\documentclass{amsart}
\usepackage{amsfonts,amssymb,epsfig}
\usepackage[latin1]{inputenc}
\usepackage[all]{xy}
\usepackage{amsmath}
\usepackage{delarray}
\usepackage{amssymb}
\newtheorem{thm}{Theorem}[section]

\newtheorem{lem}[thm]{Lemma}

\theoremstyle{definition}

\theoremstyle{remark}

\numberwithin{equation}{section}

\renewcommand{\epsilon}{\varepsilon}

\newcommand{\fg}{\mathfrak{g}}
\newcommand{\fl}{\mathfrak{l}}
\newcommand{\cF}{\mathcal F}

\newcommand{\bbR}{\mathbb R}

\newcommand{\bbN}{\mathbb N}
\newcommand{\bbZ}{\mathbb Z}

\newcommand{\N}{{\mathbb N}}
\newcommand{\Z}{{\mathbb Z}}

\newcommand{\R}{{\mathbb R}}

\newcommand{\norm}{{\|}}

\newcommand{\QED}{\hfill $\square$\vspace{2mm}}

\begin{document}

\title[Levi decomposition for smooth Poisson structures]
{Levi decomposition for smooth Poisson structures}

\author{Philippe Monnier}
\address{Departamento de Matem\'atica, Instituto Superior T\'ecnico, Lisbon, Portugal}
\email{pmonnier@math.ist.utl.pt}
\thanks{Ph. Monnier was partially supported by POCTI/FEDER and by FCT Fellowship
SFRH/BPD/12057/2003.}
\author{Nguyen Tien Zung}
\address{Laboratoire Emile Picard, UMR 5580 CNRS, UFR MIG, Université Toulouse III}
\email{tienzung@picard.ups-tlse.fr}

\keywords{Poisson structures, singular foliations, Lie algebroids,
normal forms, Levi decomposition, semisimple algebra of symmetry,
linearization}

\subjclass{53D17, 32S65}
\date{March/2004}%

\begin{abstract}
We prove the existence of a local smooth Levi decomposition for
smooth Poisson structures and Lie algebroids near a singular
point. This Levi decomposition is a kind of normal form or partial
linearization, which was established in the formal case by Wade
\cite{Wade-Levi1997} and in the analytic case by the second author
\cite{Zung-Levi2002}. In particular, in the case of smooth Poisson
structures with a compact semisimple linear part, we recover
Conn's smooth linearization theorem \cite{Conn-Smooth1985}, and in
the case of smooth Lie algebroids with a compact semisimple
isotropy Lie algebra, our Levi decomposition result gives a
positive answer to a conjecture of Weinstein
\cite{Weinstein-Linearization2000} on the smooth linearization of
such Lie algebroids. In the appendix of this paper, we show an
abstract Nash-Moser normal form theorem, which generalizes our
Levi decomposition result, and which may be helpful in the study
of other smooth normal form problems.
\end{abstract}
\maketitle

\section{Introduction}

In the study of Poisson structures, in particular their local
normal forms, one is led naturally to the problem of finding a
semisimple subalgebra of the (infinite-dimensional) Lie algebra of
functions under the Poisson bracket: such a subalgebra can be
viewed as a semisimple Lie algebra of symmetry for the
corresponding Poisson structure, and by linearizing it one get a
partial linearization of the Poisson structure, which in some case
leads to a full linearization. We call it the Levi decomposition
problem, because it is an infinite-dimensional analog of the
classical Levi decomposition for finite-dimensional Lie algebras.

Recall that, if $\mathfrak{l}$ is a finite-dimensional Lie algebra
and $\mathfrak{r}$ is the solvable radical of $\mathfrak{l}$, then
there is a semisimple subalgebra $\fg$ of $\mathfrak{l}$ such that
$\mathfrak{l}$ is a semi-direct product of $\fg$ with
$\mathfrak{r}$: $\mathfrak{l} = \fg \ltimes \mathfrak{r}$. This
semidirect product is called the Levi decomposition of
$\mathfrak{l}$, and $\fg$ is called the Levi factor of
$\mathfrak{l}$. The classical theorem of Levi and Malcev says that
$\fg$ exists and is unique up to conjugations in $\mathfrak{l}$,
see, e.g., \cite{Bourbaki-Lie1960}.

The Levi-Malcev theorem does not hold for infinite dimensional
algebras in general. But a formal version of it holds for filtered
pro-finite Lie algebras: if $\mathfrak{L} \supset \mathfrak{L}_1
\supset ... \supset \mathfrak{L}_i \supset ... $ where
$\mathfrak{L}_i$ are ideals of a Lie algebra $\mathfrak{L}$ such
that $[\mathfrak{L}_i,\mathfrak{L}_j] \subset \mathfrak{L}_{i+j}$
and $\dim \mathfrak{L}/\mathfrak{L}_i$ are finite, then the
projective limit $\lim_{i \to \infty} \mathfrak{L}/\mathfrak{L}_i$
admits a Levi factor (which is isomorphic to the Levi factor for
$\mathfrak{L}/\mathfrak{L}_1$). The proof of this formal infinite
dimensional Levi decomposition is absolutely similar to the proof
of the classical Levi-Malcev theorem. And the formal Levi
decomposition for singular foliations
\cite{Cerveau-Involutive1979} and Poisson structures
\cite{Wade-Levi1997} are instances of this infinite dimensional
formal Levi decomposition.

In \cite{Zung-Levi2002}, the second author obtained the local
analytic Levi decomposition theorem for analytic Poisson
structures which vanish at a point. This theorem generalizes
Conn's linearization theorem for analytic Poisson structure with a
semisimple linear part \cite{Conn-Analytic1984}, and is at the
base of some new analytic linearization results for Poisson
structures and Lie algebroids \cite{Zung-Levi2002,DuZu-Affn2002}.

The aim of this paper is to establish the local smooth Levi
decomposition theorem for smooth Poisson structures and Lie
algebroids which vanish at a point. Our main theorem (Theorem
\ref{thm:LeviPoissonS}) is a generalization of Conn's smooth
linearization theorem \cite{Conn-Smooth1985} for Poisson
structures with a compact semisimple linear part, and provides a
local smooth semi-linearization for any smooth Poisson structure
whose linear part (when considered as a Lie algebra) contains a
compact semisimple subalgebra.

Let $\Pi$ be a $C^{p}$ Poisson structure ($p \in \bbN \cup
\{\infty\}$) in a neighborhood of $0$ in $\bbR^n$, which vanishes
at the origin. Denote by $\mathfrak{l}$ the $n$-dimensional Lie
algebra of linear functions in $\bbR^n$ under the Lie-Poisson
bracket $\Pi_1$ which is the linear part of $\Pi$ at 0, and by
$\mathfrak{g}$ a compact semisimple subalgebra of $\mathfrak{l}$.
(Without loss of generality one can assume that $\mathfrak{g}$ is
a maximal compact semisimple subalgebra of $\mathfrak{l}$, and we
will call $\mathfrak{g}$ a {\it compact Levi factor} of
$\mathfrak{l}$). Denote by $(x_1,...,x_m,y_1,...,y_{n-m})$ a
linear basis of $\mathfrak{l}$, such that $x_1,...,x_m$ span $\fg$
($\dim \fg = m$), and $y_1,...,y_{n-m}$ span a linear complement
$\mathfrak{r}$ of $\fg$ in $\mathfrak{l}$ which is invariant under
the adjoint action of $\fg$. Denote by $c_{ij}^k$ and $a_{ij}^k$
the structural constants of $\fg$ and of the action of $\fg$ on
$\mathfrak{r}$ respectively: $[x_i,x_j] = \sum_k c_{ij}^k x_k$ and
$[x_i,y_j] = \sum_k a_{ij}^k y_k$. We say that $\Pi$ admits a
local \emph{$C^q$-smooth Levi decomposition with respect to $\fg$}
if there exists a local $C^q$-smooth system of coordinates
$(x^{\infty}_1,...,x^{\infty}_m,y^{\infty}_1,...,
y^{\infty}_{n-m})$, with $x^{\infty}_{i} = x_i +$ higher order
terms and $y^{\infty}_{i} = y_i +$ higher order terms, such that
in this coordinate system the Poisson structure has the form
\begin{equation}
\label{eqn:LeviPoissonS1} \Pi = \frac{1}{2} \big[\sum c_{ij}^k
x^{\infty}_{k} \frac{\partial}{\partial x^{\infty}_{i}} \wedge
\frac{\partial }{ \partial x^{\infty}_{j}} + \sum a_{ij}^k
y^{\infty}_{k} \frac{\partial }{
\partial x^{\infty}_{i}} \wedge \frac{\partial }{ \partial
y^{\infty}_{j}} + \sum F_{ij} \frac{\partial }{
\partial y^{\infty}_{i}} \wedge \frac{\partial }{ \partial y^{\infty}_{j}} \big]
\end{equation}
where $F_{ij}$ are some functions in a neighborhood of $0$ in
$\bbR^n$. In other words, we have
\begin{equation}
\label{eqn:LeviPoissonS2} \{x^{\infty}_{i},x^{\infty}_{j}\} = \sum
c_{ij}^k x^{\infty}_{k} \ {\rm and} \
\{x^{\infty}_{i},y^{\infty}_{j}\} = \sum a_{ij}^k y^{\infty}_{k},
\end{equation}
i.e. the functions $x^\infty_1,\dots,x^\infty_m$ span a compact
Levi factor (isomorphic to $\mathfrak{g}$) and their Hamiltonian
vector fields $X_{x^\infty_1},\dots, X_{x^\infty_m}$ are linear in
the coordinate system
$(x^{\infty}_1,...,x^{\infty}_m,y^{\infty}_1,...,
y^{\infty}_{n-m})$.

\begin{thm}
\label{thm:LeviPoissonS} There exists a positive integer $l$
(which depends only on the dimension $n$) such that any
$C^{2q-1}$-smooth Poisson structure $\Pi$ in a neighborhood of $0$
in $\bbR^n$ which vanishes at $0$, where $q \in \bbN \cup
\{\infty\}$, $q \geq l$, admits a local $C^{q}$-smooth Levi
decomposition (with respect to any compact semisimple Lie
subalgebra $\fg$ of the Lie algebra $\mathfrak{l}$ which
corresponds to the linear part of $\Pi$ at $0$).
\end{thm}

A particular case of the above theorem is when $\fg = \fl$, i.e.
when the linear part of $\Pi$ is compact semisimple. In this case
a local Levi decomposition is nothing but a local linearization of
the Poisson structure, and we recover the smooth linearization
theorem of Conn \cite{Conn-Smooth1985} for a smooth Poisson
structure with a compact semisimple linear part. When $\fl = \fg
\oplus \bbR$, a Levi decomposition is still a linearization of
$\Pi$. In general, one may consider a Levi decomposition (we also
call it a \emph{Levi normal form}, see \cite{Zung-Levi2002}) as a
partial linearization of $\Pi$.

Similarly to the analytic case \cite{Zung-Levi2002}, an analogue
of Theorem \ref{thm:LeviPoissonS} holds for smooth Lie algebroids:

\begin{thm}
\label{thm:LeviAlgebroidS} Let $A$ be a local $N$-dimensional
$C^{2q-1}$-smooth Lie algebroid over $(\bbR^n,0)$ with the anchor
map $\# : A \to T\bbR^n$, such that $\#a = 0$ for any $a \in A_0$,
the fiber of $A$ over point $0$, where $q=\infty$ or is a natural
number which is large enough ($q \geq l$, where $l$ is a natural
number which depends only on $N$ and $n$). Denote by
$\mathfrak{l}$ the $N$-dimensional Lie algebra in the linear part
of $A$ at $0$ (i.e. the isotropy algebra of $A$ at $0$), and by
$\fg$ a compact semisimple Lie subalgebra of $\mathfrak{l}$. Then
there exists a local $C^q$-smooth system of coordinates
$(x^\infty_1,...,x^\infty_n)$ of $(\bbR^n,0)$, and a local
$C^q$-smooth basis of sections
$(s^\infty_1,s^\infty_2,...,s^\infty_m,v^\infty_1,...,v^\infty_{N-m})$
of $A$, where $m = \dim \fg$, such that we have:
\begin{equation}
\begin{array}{l}
[s^\infty_i,s^\infty_j] = \sum_k c_{ij}^k s^\infty_k  \ , \cr
[s^\infty_i, v^\infty_j] = \sum_k a_{ij}^k v^\infty_k \ ,\cr \#
s^\infty_i = \sum_{j,k} b_{ij}^k x^\infty_k \partial /
\partial x^\infty_j \ ,
\end{array}
\end{equation}
where $c_{ij}^k,a_{ij}^k,b_{ij}^k$ are constants, with $c_{ij}^k$
being the structural constants of the compact semisimple Lie
algebra $\fg$.
\end{thm}

The meaning of the above theorem is that the algebra of sections
of $A$ admits a Levi factor (Lie isomorphic to $\fg$), spanned by
$s^\infty_1,s^\infty_2,...,s^\infty_m$, whose action can be
linearized. Theorem \ref{thm:LeviAlgebroidS} is called the local
smooth Levi decomposition theorem for smooth Lie algebroids. As a
particular case of this theorem, we obtain the following result,
conjectured by A. Weinstein \cite{Weinstein-Linearization2000}:
any smooth Lie algebroid whose anchor vanishes at a point and
whose corresponding isotropy Lie algebra at that point is compact
semisimple is locally smoothly linearizable.

Remark that, compared to the analytic case, in the smooth case
considered in \cite{Conn-Smooth1985} and in the present paper we
need the additional condition of compactness on our semisimple Lie
(sub)algebra $\fg$. In a sense, this compactness condition is
necessary, due to the following result of Weinstein
\cite{Weinstein-Real1987}: any real semisimple Lie algebra of real
rank at least 2 is smoothly degenerate, i.e. there is a smoothly
nonlinearizable Poisson structure with a linear part corresponding
to it.

We hope that the results of this paper will be useful for finding
new smoothly nondegenerate Lie algebras (and Lie algebroids) in
the sense of Weinstein \cite{Weinstein-Poisson1983}. In
particular, our smooth Levi decomposition is one of the main steps
in the study of smooth linearizability of Poisson structures whose
linear part corresponds to a real semisimple Lie algebra of real
rank 1 (this case was left out by Weinstein
\cite{Weinstein-Real1987}). This problem will be studied in a
separate work.

Our proof of Theorem \ref{thm:LeviPoissonS} is based on the
Nash-Moser fast convergence method (see, e.g.,
\cite{Hamilton-NashMoser1982}) applied to Fréchet spaces of smooth
functions and vector fields. In particular, our algorithm for
constructing a convergent sequence of smooth coordinate
transformations, which is a combination of smoothing operators
with the algorithm in \cite{Zung-Levi2002} for the analytic case,
is inspired by Hamilton's ``near projections'' in his proof of the
so-called Nash-Moser theorem for exact sequences
\cite{Hamilton-Complex1977}. Besides smoothing operators for tame
Fr\'echet spaces, we will need homotopy operators for certain
Chevalley-Eilenberg complexes with vanishing first and second
cohomologies. The homotopy operators and the smoothing operators
are both already present in Conn's paper \cite{Conn-Smooth1985},
and in a sense the present paper is a further development of
\cite{Conn-Smooth1985} and follows more or less the same
organization.

Using the fact that Lie algebroids can be viewed as fiber-wise
linear Poisson structures, one can immediately deduce Theorem
\ref{thm:LeviAlgebroidS} from the proof given below of Theorem
\ref{thm:LeviPoissonS}, simply by restricting some functional
spaces, in a way absolutely similar to the analytic case (see
Section 6 of \cite{Zung-Levi2002}). That's why we will mention
only briefly the proof of Theorem \ref{thm:LeviAlgebroidS}, after
the full proof Theorem \ref{thm:LeviPoissonS}.

The rest of this paper, except the appendix, is devoted mainly to
the proof of Theorem \ref{thm:LeviPoissonS}, and is organized as
follows. In Sections \ref{section:homotopy} and
\ref{section:inequalities} we write down important inequalities
involving homotopy operators and smoothing operators that will be
used. Then in Section \ref{section:proof} we present our algorithm
for constructing the required new systems of coordinates, and give
a proof of Theorem \ref{thm:LeviPoissonS}, modulo some technical
lemmas. These lemmas are proved in Section
\ref{section:technical}. In Section \ref{section:Algebroids} we
briefly explain how to modify (in an obvious way) the proof of
Theorem \ref{thm:LeviPoissonS} to get a proof of Theorem
\ref{thm:LeviAlgebroidS}.

In the appendix, we present an abstract Nash-Moser smooth normal
form theorem, which generalizes Theorems \ref{thm:LeviPoissonS}
and \ref{thm:LeviAlgebroidS}. We hope that this abstract normal
form theorem can be used or easily adapted for the study of other
smooth normal form problems (of functions, dynamical systems,
various geometric structures, etc.).

\section{Homotopy operators}
\label{section:homotopy}

Similarly to the analytic case
\cite{Conn-Analytic1984,Zung-Levi2002}, in order to prove Theorem
\ref{thm:LeviPoissonS}, we will need a normed version of
Whitehead's lemma about the vanishing of cohomology of the
semisimple algebra $\fg$, with respect to certain orthogonal
modules of $\fg$ constructed below. Our modules will be spaces of
real functions or vector fields, equipped with Sobolev norms, and
the action of $\fg$ will preserve these norms.

Consider a Lie algebra $\mathfrak{l}$ of dimension $n$ together
with a compact semisimple Lie subalgebra $\mathfrak{g} \subset
\mathfrak{l}$ of dimension $m$. (Our Poisson structure will live
in a neighborhood of $0$ in the dual space $\bbR^n =
\mathfrak{l}^*$ of $\mathfrak{l}$). Denote by $G$ the
simply-connected compact semisimple Lie group whose Lie algebra is
$\mathfrak{g}$. Then $G$ acts on $\bbR^n = \mathfrak{l}^*$ by the
coadjoint action. Since $G$ is compact, we can fix a linear
coordinate system $(x_1,\dots,x_m,y_1,\dots,y_{n-m})$ such that
the Euclidean metric on $\bbR^n$ with respect to this coordinate
system is invariant under the action of $G$, and the first $m$
coordinate $(x_1,\dots,x_m)$ come from $\mathfrak{g}$. In other
words, there is a basis $(\xi_1,\dots,\xi_m)$ of $\mathfrak{g}$
such that each $\xi_i$, considered as an element of $\mathfrak{l}$
and viewed as a linear function on $\mathfrak{l}^*$, gives rise to
the coordinate $x_i$.



For each positive number $r > 0$, denote by $B_r$ the closed ball
of radius $r$ in $\R^n$ centered at $0$. The group $G$ (and hence
the algebra $\fg$) acts linearly on the space of functions on
$B_r$ via its action on $B_r$: for each function $F$ and element
$g \in G$ we put
\begin{equation}
\label{eqn:G-action-function} g(F) (z) := F (g^{-1}(z)) = F
(Ad^*_{g^{-1}}z) \ \ \forall \ z \in B_r .
\end{equation}

For each nonnegative integer $k \geq 0$ and each pair of
real-valued functions $F_1,F_2$ on $B_r$, we will define the
Sobolev inner product of $F_1$ with $F_2$ with respect to the
Sobolev $H_k$-norm as follows:

\begin{equation}
\label{eqn:H-inner} \langle F_1,F_2\rangle^H_{k,r} :=
\sum_{|\alpha| \leq k} \int_{B_r} \left(\frac{|\alpha|!}{\alpha
!}\right) \left( \frac{\partial^{|\alpha|}F_1}{\partial z^\alpha }
(z) \right)  \left( \frac{\partial^{|\alpha|}F_2}{\partial
z^\alpha } (z) \right) d\mu(z) ,
\end{equation}
where $d\mu$ is the standard Lebesgue measure on $\R^n$. The
Sobolev $H_k$-norm of a function $F$ on $B_r$ is
\begin{equation}
\|F\|^H_{k,r} := \sqrt{\langle F,F\rangle_{k,r}} \ .
\end{equation}

We will denote by $\mathcal{C}_r$ the subspace of the space
$C^\infty(B_r)$ of $C^\infty$-smooth real-valued functions on
$B_r$, which consists of functions vanishing at $0$ whose first
derivatives also vanish at $0$. Then the action of $G$ on
$\mathcal{C}_r$ defined by (\ref{eqn:G-action-function}) preserves
the Sobolev inner products (\ref{eqn:H-inner}).

Denote by $\mathcal{Y}_r$ the space of $C^\infty$-smooth vector
fields on $B_r$ of the type
\begin{equation}
u = \sum_{i=1}^{n-m} u_i \partial / \partial y_i \ ,
\end{equation}
such that $u_i$ vanish at $0$ and their first derivatives also
vanish at $0$.

Recall that $(\xi_1,\hdots,\xi_m)$ is the basis of $\fg$ which
correspond to the coordinates $(x_1,\hdots,x_m)$ on $\bbR^n =
\mathfrak{l}^*$. The space $\mathcal{Y}_r$ is a $\fg$-module under
the following action:
\begin{equation}
\xi_i \cdot \sum_j u_j \partial / \partial y_j := \Big[ \sum_{jk}
c_{ij}^k x_k \frac{\partial }{
\partial x_j} + \sum_{jk} a_{ij}^k y_k \frac{\partial }{
\partial y_j}\ ,\ \sum_j u_j \partial /
\partial y_j \Big] \ ,
\end{equation}
where $X_i = \sum_{jk} c_{ij}^k x_k \partial / \partial x_j +
\sum_{jk} a_{ij}^k y_k
\partial / \partial y_j$ are the linear vector fields which generate the linear
orthogonal coadjoint action of $\fg$ on $\bbR^n$.

Equip $\mathcal{Y}_r$ with Sobolev inner products:
\begin{equation}
\langle u,v\rangle^H_{k,r} := \sum_{i=1}^{n-m} \langle u_i,v_i
\rangle_{k,r} \ ,
\end{equation}
and denote by $\mathcal{Y}^H_{k,r}$ the completion of
$\mathcal{Y}_r$ with respect to the corresponding $H_{k,r}$-norm.
Then $\mathcal{Y}^H_{k,r}$ is a separable real Hilbert space on
which $\fg$ and $G$ act orthogonally.

The following infinite dimensional normed version of Whitehead's
lemma is taken from Proposition 2.1  of \cite{Conn-Smooth1985}:

\begin{lem}[Conn]
\label{lem:HomotopyOp} For any given positive number $r$, and $W =
\mathcal{C}_r$ or $\mathcal{Y}_r$ with the above action of $\fg$,
consider the (truncated) Chevalley-Eilenberg complex
$$W \stackrel{\delta_0}{\rightarrow} W \otimes \wedge^1 \fg^* \stackrel{\delta_1}{\rightarrow}
W \otimes \wedge^2 \fg^* \stackrel{\delta_2}{\rightarrow} W
\otimes \wedge^3 \fg^* .
$$ Then there is a chain of operators
$$W \stackrel{h_0}{\leftarrow} W \otimes \wedge^1 \fg^* \stackrel{h_1}{\leftarrow} W \otimes
\wedge^2 \fg^* \stackrel{h_2}{\leftarrow} W \otimes \wedge^3 \fg^*
$$ such that
\begin{equation}
\label{eqn:homotopy}
\begin{array}{c}
\delta_0 \circ h_0 + h_1 \circ \delta_1 = {\rm Id}_{W \otimes \wedge^1\fg^*} \ , \\
\delta_1 \circ h_1 + h_2 \circ \delta_2 = {\rm Id}_{W \otimes
\wedge^2\fg^*} \ .
\end{array}
\end{equation}
Moreover, there exist a constant $C > 0$, which is independent of
the radius $r$ of $B_r$, such that
\begin{equation}
\label{eqn:estimate-h2} \|h_j(u)\|^H_{k,r} \leq C \|u\|^H_{k,r} \,
, \ \ j=0,1,2
\end{equation}
for all $k \geq 0$ and $u \in W \otimes \wedge^{j+1} \fg^*$. If
$u$ vanishes to an order $l \geq 0$ at the origin, then so does
$h_j(u)$.
\end{lem}

{\it Proof.} Strictly speaking, Conn \cite{Conn-Smooth1985} only
proved the above lemma in the case when $\fg = \mathfrak{l}$ and
for the module $\mathcal{C}_r$, but his proof is quite general and
works perfectly in our situation without any modification. Here,
we will just recall the main idea of this proof. The action of
$\fg$ on $W$ can be extended to the completion ${\tilde W}_k$ of
$W$ with respect to the $H_{k,r}$-norm (this is the Sobolev space
$H_k(B_r)$ when $W = \mathcal{C}_r$ and $\mathcal{Y}_{k,r}^H$ when
$W = \mathcal{Y}_r$). We can decompose ${\tilde W}_k$ as an
orthogonal direct sum of $\fg$-modules ${\tilde
W}_k^0\oplus{\tilde W}_k^1$ where ${\tilde W}_k^0$ is a trivial
$\fg$-module and ${\tilde W}_k^1$ can be decomposed as a Hilbert
direct sum of finite dimensional irreducible $\fg$-invariant
subspaces. This decomposition induces a decomposition $W=W^0\oplus
W^1$. We can construct a homotopy operator $h_i^\prime :
W^0\otimes \wedge^{i+1}\fg\longrightarrow W^0\otimes
\wedge^{i}\fg$ by tensoring the identity mapping of $W^0$ with a
homotopy operator for the trivial $\fg$-module $\bbR$. To
construct the homotopy operator $h_i^{\prime\prime}$ on
$W^1\otimes \wedge^{i+1}\fg$, we can restrict to the case when
$W^1$ is irreducible. Then we define the $h_i^{\prime\prime}$ by
\begin{eqnarray*}
h_0^{\prime\prime}(w) &=& \Gamma^{-1}\cdot (\sum_i \xi_k\cdot w(\xi_k))\\
h_1^{\prime\prime}(w) &=& \sum_i\xi^*\otimes(\Gamma^{-1}\cdot (\sum_k \xi_k\cdot w(\xi_i\wedge\xi_k)))\\
h_2^{\prime\prime}(w) &=& \sum_{ij}\xi_i^*\wedge\xi_j^* \otimes (
\Gamma^{-1}\cdot (\sum_k \xi_k\cdot
w(\xi_i\wedge\xi_j\wedge\xi_k)))
\end{eqnarray*}
where $\{\xi_i^*\}$ is the dual basis of $\{\xi_i\}$ and $\Gamma$
is the Casimir element of $\fg$. Then one can show that
$$
\|h_i^{\prime\prime}(w)\|_{k,r}^H\leq C \|w\|_{k,r}^H
$$
with $C=m(\min_{\gamma\in\mathcal{J}}\|\gamma\|)^{-1}$, where
$\mathcal{J}$ is the weight lattice of $\fg$. \QED

For simplicity, in the sequel we will denote the homotopy
operators $h_j$ in the above lemma simply by $h$. Relation
(\ref{eqn:homotopy}) will be rewritten simply as follows:
\begin{equation}
\label{eqn:homotopy'} {\rm Id} - \delta \circ h =  h \circ \delta
\ .
\end{equation}
The meaning of the last equality is as follows: if $u$ is an
1-cocycle or 2-cocycle, then it is also a coboundary, and $h(u)$
is an explicit primitive of $u$: $\delta(h(u)) = u$. If $u$ is a
``near cocycle'' then $h(u)$ is also a ``near primitive'' for $u$.

For convenience, in the sequel, instead of Sobolev norms, we will
use the following absolute forms:
\begin{equation}
\label{eqn:absolute-norm1} \|F\|_{k,r} := \sup_{|\alpha| \leq k}
\sup_{z \in B_r} |D^\alpha F(z)|
\end{equation}
for $F \in \mathcal{C}_r$, where the sup runs over all partial
derivatives of degree $|\alpha|$ at most $k$. More generaly, if
$F=(F_1,\hdots,F_m)$ is a smooth maping from $B_r$ to $\bbR^m$ we
can define
\begin{equation}
\label{eqn:absolute-norm3} \|F\|_{k,r} := \sup_i \sup_{|\alpha|
\leq k} \sup_{z \in B_r} |D^\alpha F_i(z)|\,.
\end{equation}
Similarly, for $u = \sum_{i=1}^{n-m} u_i \partial / \partial y_i
\in \mathcal{Y}_r$ we put
\begin{equation}
\label{eqn:absolute-norm2} \|u\|_{k,r} := \sup_{i} \sup_{|\alpha|
\leq k} \sup_{z \in B_r} |D^\alpha u_i(z)|\,.
\end{equation}
The absolute norms $\|.\|_{k,r}$ are related to the Sobolev norms
$\|.\|^H_{k,r}$ as follows:
\begin{equation}
\label{eqn:2norms} \|F\|_{k,r} \leq C_1 \|F\|^H_{k+s,r} \ {\rm
and} \ \|F\|^H_{k,r} \leq C_2 (n+1)^k  \|F\|_{k,r}
\end{equation}
for any $F$ in $\mathcal{C}_r$ or $\mathcal{Y}_r$ and any $k \geq
0$, where $s = [\frac{n}{2}] + 1$ and $C_1$ and $C_2$ are positive
constants which do not depend on $k$. A priori, the constants
$C_1$ and $C_2$ depend continuously on $r$ (and on the dimension
$n$), but later on we will always assume that $1 \leq r \leq 2$,
and so may assume $C_1$ and $C_2$ to be independent of $r$. The
above first inequality is a version of the classical Sobolev's
lemma for Sobolev spaces.
The second inequality follows directly from the definitions of the
norms. Combining it with Inequality (\ref{eqn:estimate-h2}), we
obtain the following estimate for the homotopy operators $h$ with
respect to absolute norms:

\begin{equation}
\label{eqn:estimate-h1} \|h(u)\|_{k,r} \leq C (n+1)^{k+s}
\|u\|_{k+s,r}
\end{equation}
for all $k \geq 0$ and $u \in W \otimes \wedge^{j+1} \fg^*$
($j=0,1,2$), where $W = \mathcal{C}_r$ or $\mathcal{Y}_r$. Here $s
= [\frac{n}{2}] + 1$, $C$ is a positive constant which does not
depend on $k$ (and on $r$, provided that $1 \leq r \leq 2$).

\section{Smoothing operators and some useful inequalities}
\label{section:inequalities}

We will refer to \cite{Hamilton-NashMoser1982} for the theory of
tame Fr\'echet spaces used here. It is well-known that the space
$C^\infty(B_r)$ with absolute norms (\ref{eqn:absolute-norm1}) is
a tame Fr\'echet space. Since $\mathcal{C}_r$ is a tame direct
summand of $C^\infty(B_r)$, it is also a tame Fr\'echet space.
Similarly, $\mathcal{Y}_r$ with absolute norms
(\ref{eqn:absolute-norm2}) is a tame Fr\'echet space as well. In
particular, $\mathcal{C}_{r}$ and $\mathcal{Y}_{r}$  admit
smoothing operators and interpolation inequalities:

For each $t > 1$ there is a linear operator $S(t) = S_r(t)$ from
$\mathcal{C}_{r}$ to itself, with the following properties:
\begin{equation}
\label{eqn:smoothing1} \|S(t) F\|_{p,r} \leq C_{p,q} t^{(p-q)} \|
F\|_{q,r}
\end{equation}
and
\begin{equation}
\label{eqn:smoothing2} \|(I - S(t)) F\|_{q,r} \leq C_{p,q}
t^{(q-p)} \| F\|_{p,r}
\end{equation}
for any $F \in \mathcal{C}_r$, where $p,q$ are any nonnegative
integers such that $p \geq q$, $I$ denotes the identity map, and
$C_{p,q}$ denotes a constant which depends on $p$ and $q$.

The second inequality means that $S(t)$ is close to identity and
tends to identity when $t \to \infty$. The first inequality means
that $F$ becomes ``smoother'' when we apply $S(t)$ to it. For
these reasons, $S(t)$ is called the smoothing operator.

{\it Remark.} Some authors write $e^{t(p-q)}$ and $e^{t(q-p)}$
instead of $t^{(p-q)}$ and $t^{(p-q)}$ in the above inequalities.
The two conventions are related by a simple rescaling $t=e^\tau.$

There is a similar smoothing operator from $\mathcal{Y}_r$ to
itself, which by abuse of language we will also denote by $S(t)$
or $S_r(t)$. We will assume that inequalities
(\ref{eqn:smoothing1}) and (\ref{eqn:smoothing2}) are still
satisfied when $F$ is replaced by an element of $\mathcal{Y}_r$.

For any $F$ in $\mathcal{C}_r$ or $\mathcal{Y}_r$, and nonnegative
integers $p_1 \geq p_2 \geq p_3$, we have the following
interpolation estimate:
\begin{equation}
\label{eqn:interpolation} (\|F\|_{p_2,r})^{p_3-p_1} \leq
C_{p_1,p_2,p_3} (\|F\|_{p_1,r})^{p_3-p_2} (\|F\|_{p_3,r})^{p_2 -
p_1}
\end{equation}
where $C_{p_1,p_2,p_3}$ is a positive constant which may depend on
$p_1,p_2,p_3$.

 {\it Remark.} A priori, the constants $C_{p,q}$ and $C_{p_1,p_2,p_3}$
also depend on the radius $r$. But later on, we will always have
$1 \leq r \leq 2$ and so we may choose them to be independent of
$r$.\\

In the proof of Theorem \ref{thm:LeviPoissonS}, we will use local
diffeomorphisms of $\R^n$ of type $Id+\chi$ where $\chi(0)=0$, and
$Id$ denotes the identity map from $\bbR^n$ to itself. The
following lemmas allow to control operations on this kind of
diffeomorphisms as the composition with a map or the inverse.

\begin{lem}
Let $r$ and $\eta<1$ be two strictly positive real numbers.
Consider a smooth map $\Phi : B_r\rightarrow \bbR^n$ of the type
$Id+\chi$ with $\chi(0)=0$. Suppose that $\|\chi\|_{1,r}<\eta$.
Then we have
\begin{equation}
B_{r(1-\eta)}\subset \Phi(B_r) \subset B_{r(1+\eta)}\,.
\end{equation}
\label{lem:boules}
\end{lem}

{\it Proof : } According to the hypotheses we have $\|\chi(x)\| <
\eta\|x\|$ for every $x$ in $B_r$. Therefore, we can write
$\|\Phi(x)\| < (1+\eta)r$ and so, $\Phi(B_r)\subset
B_{r(1+\eta)}$.

Now, we consider the map ${\hat \Phi} : B_{r(1+\eta)}\rightarrow
B_{r(1+\eta)}$ which is $\Phi$ on $B_r$ and is defined on
$B_{r(1+\eta)}\setminus B_r$ as
follows.\\
Let $x$ be such that $\|x\|=r$. We consider
$x_1=\frac{2+\eta}{2}x$ and $x_2=(1+\eta)x$. If $z=\lambda
x+(1-\lambda)x_1$ with $0\leq \lambda\leq 1$, then ${\hat
\Phi}(z)=\lambda\Phi(x)+(1-\lambda)x$. If $z=\lambda
x_1+(1-\lambda)x_2$ then ${\hat \Phi}(z)=\lambda
x+(1-\lambda)x_2$.

 This map is continous and is the identity on the boundary of $B_{r(1+\eta)}$.
According to Brouwer's theorem, the image of ${\hat \Phi}$ is
$B_{r(1+\eta)}$.

Now, note that if $z=\lambda x+(1-\lambda)x_1$ with
$0\leq\lambda\leq 1$ then we have $$\|{\hat
\Phi}(z)\|=\|x+\lambda\chi(x)\|\geq \|x\|-\lambda\|\chi(x)\|\,.$$
Therefore, $\|{\hat \Phi}(z)\|>r(1-\eta)$.

Moreover, if $v=\lambda x_1+(1-\lambda)x_2$ with $0\leq\lambda\leq
1$ then we have
$$\|{\hat \Phi}(v)\|=\|\lambda x+(1-\lambda)(1+\eta)x\|=r(1+\eta(1-\lambda))>r\,.$$

We deduce that if $z=\lambda x+(1-\lambda)x_1$ with
$0\leq\lambda\leq 1$ then we have $y$ is in $B_{r(1-\eta)}$ then,
we have $y={\hat \Phi}(z)$ ($z$ is, a priori, in $B_{r(1+\eta)}$)
and $z$ must be in $B_r$. Consequently, $y$ is in the image of
$\Phi$. \QED

\begin{lem}[\cite{Conn-Smooth1985}]
\label{lem:composition} Let $r > 0$ and $1 > \eta > 0$ be two
positive numbers. Consider two smooth maps
$$f:B_{r(1+\eta)}\rightarrow\bbR^q\quad {\mbox { and }}\quad \chi:B_r\rightarrow \bbR^n$$
(where the closed balls $B_r$ and $B_{r(1+\eta)}$ are in $\bbR^n$,
and $q$ is a natural number) such that $\chi(0)=0$ and
$\|\chi\|_{1,r}<\eta$. Then the composition $f\circ (id+\chi)$ is
a smooth map from $B_r$ to $\bbR^n$ which satisfies the following
inequalities:
\begin{eqnarray}
\|f\circ(id+\chi)\|_{k,r} &\leq& \|f\|_{k,r(1+\eta)}(1+P_k(\|\chi\|_{k,r}))\\
\|f\circ(id+\chi)-f\|_{k,r} &\leq&
Q_k(\|\chi\|_{k,r})\|f\|_{k,r(1+\eta)}+ M\|\chi\|_{0,r}
\|f\|_{k+1,r(1+\eta)}
\end{eqnarray}
where $M$ is a positive constant and $P_k(t)$,$Q_k(t)$ are
polynomials of degree $k$ with vanishing constant term (and which
are independent of $f$ and $\chi$).
\end{lem}

The proof of the above lemma, which can be found in
\cite{Conn-Smooth1985}, is straightforward and is based solely on
the Leibniz rule of derivation. We will call inequalities such as
in the above lemma \emph{Leibniz-type inequalities}. Similarly, we
have another Leibniz-type inequality, given in the following
lemma.

\begin{lem}
\label{lem:compositionL} With the same hypotheses as in the
previous lemma, we have
\begin{eqnarray}\label{eqn:compL}
\|f\circ (id+\chi)\|_{2k-1,r} &\leq&  \|f\|_{2k-1,r(1+\eta)} P_k(\|\chi\|_{k,r}) \\
  &+& \|\chi\|_{2k-1,r}  \|f\|_{k,r(1+\eta)}  Q_k(\|\chi\|_{k,r})
\,,\nonumber
\end{eqnarray}
where $P_k(t)$ and $Q_k(t)$ are polynomials (which are independent
of $f$ and $\chi$).
\end{lem}

{\it Proof.} Denote by $\theta$ the map $Id+\chi$. If $I$ is a
multiindex such that $|I|\leq 2k-1$ ($|I|$ denotes the sum of the
components of $I$), it is easy to show, by induction on $|I|$,
that
$$\frac{\partial^{|I|} (f\circ\theta)}{\partial
x^{|I|}}=\sum_{1\leq|\alpha|\leq|I|}\big(\frac{\partial^{|\alpha|}
f}{\partial x^\alpha}\circ\theta\big)A_\alpha(\theta)\,,$$ where
$A_\alpha(\theta)$ is of the type
\begin{equation}\label{eqn:Aalpha}
A_\alpha(\theta)= \sum_{\begin{array}{cc}
{\scriptstyle 1\leq u_i\leq n\,,\, |\beta_i|\geq 1}\\
{\scriptstyle |\beta_1|+\hdots+|\beta_{|\alpha|}|=|I|}
\end{array}}
a_{\beta u}\frac{\partial^{|\beta_1|}\theta_{u_1}}{\partial
x^{\beta_1}}\hdots\frac{\partial^{|\beta_{|\alpha|}|}\theta_{u_{\alpha}}}{\partial
x^{\beta_{|\alpha|}}}
\end{equation}
where $\theta_{u_1}$ is the $u_1$-component of $\theta$ and the
$a_{\beta u}$ are nonnegative integers.


We may write
$$\frac{\partial^{|I|} (f\circ\theta)}{\partial
x^{|I|}}=\sum_{k<|\alpha|\leq|I|}\big(\frac{\partial^{|\alpha|}
f}{\partial x^\alpha}\circ\theta\big)A_\alpha(\theta)+
\sum_{1\leq|\alpha|\leq k }\big(\frac{\partial^{|\alpha|}
f}{\partial x^\alpha}\circ\theta\big)A_\alpha(\theta)\,.$$ When
$k<|\alpha|\leq |I| \leq 2k -1$, all the $|\beta_i|$ in the sum
(\ref{eqn:Aalpha}) defining $A_\alpha(\theta)$ are smaller than
$k$. This gives the first term of the right hand side of
Inequality (\ref{eqn:compL}). On the other hand, when $1\leq
|\alpha| \leq k$, then in each product in the expression
(\ref{eqn:Aalpha}) of $A_\alpha(\theta)$ there is at most one
factor $\frac{\partial^{|\beta|}\theta_{u}}{\partial x^{\beta}}$
with $|\beta|>k$ (the others have $|\beta| \leq k$. This gives the
second term of the right hand side term of inequality
(\ref{eqn:compL}), and the lemma follows. \QED

\begin{lem}
Let $r>0$ be a real number and $k\geq 1$ a positive integer. There
exist a positive real number $\eta<1$ and a polynomial $P_k(t)$
such that if $\Phi: B_r \rightarrow \bbR^n$ is a smooth map of the
type $Id+\chi$ with $\chi(0)=0$ and $\|\chi\|_{0,r} < \eta$ then
$\Phi$ is a smooth local diffeomorphism which possesses an inverse
$\Psi = \Phi^{-1}$ of the type $Id+\xi$ with $\xi(0)=0$, which is
defined on (a set containing) $B_{r(1-\eta)}$ and satisfies the
following inequality:
\begin{equation}
\|\xi\|_{2k-1,r(1-\eta)}\leq
\|\chi\|_{2k-1,r}P_k(\|\chi\|_{k,r})\,.
\label{eqn:estimation-inverse}
\end{equation}
\label{lem:inverse}
\end{lem}

{\it Proof.} We choose the constant $\eta$ such that for every
smooth map $Id+\chi : B_r \rightarrow \bbR^n$ such that
$\|\chi\|_{1,r} < \eta$, the Jacobian matrix of $Id+\chi$ is
invertible at each point of $B_r$.

If $\Phi$ is a smooth map as in the theorem, according to the
inverse function theorem, it is a local diffeomorphism and has an
inverse $\Psi=Id+\xi$ which is smooth on $B_{r(1+\eta)}$ (see
Lemma \ref{lem:boules}).

Since $\Phi\circ\Psi=Id$, denoting $\Psi=(\Psi_1,\hdots,\Psi_n)$
(and the same thing for $\Phi$), we can write
$$
\frac{\partial \Psi_i}{\partial x_j}=\frac{Pol(\{\frac{\partial
\Phi_u}{\partial x_v}\})}{Jac\,\Phi}\circ \Psi
$$
where $Jac\,\Phi$ is the Jacobian determinant of $\Phi$ and
$Pol(\{\frac{\partial \Phi_u}{\partial x_v}\})$ is a homogeneous
polynomial in the $\{\frac{\partial \Phi_u}{\partial x_v}\}_{uv}$
of degree $n-1$.

By induction, we can see that for all $\alpha\in\Z_+^n$ with
$|\alpha|=\sum \alpha_i>0$ we can write (trying to simplify the
writing)
$$
\frac{\partial^{|\alpha|} \Psi_i}{\partial x^\alpha}=
\sum_{\begin{array}{cc}
{\scriptstyle 1\leq |\beta_l|\leq |\alpha|\,,\, p\leq |\alpha|+1}\\
{\scriptstyle \sum_l (|\beta_l|-1)=|\alpha|-1}
\end{array}}
\Big[\frac{a_{\beta,p}}{(Jac
\,\Phi)^p}\times\frac{\partial^{|\beta_1|} \Phi_{u_1}}{\partial
x^{\beta_1}}\hdots\frac{\partial^{|\beta_k|} \Phi_{u_k}}{\partial
x^{\beta_k}}\Big]\circ\Psi
$$
where the $a_{\beta,p}$ are non negative integers. In this
formula, the term $Jac\,\Phi$ is bounded on $B_r$, for instance,
$0<b\leq |Jac\,\Phi(z)|\leq c <1$ for all $z$ in $B_r$. This
formula is not very explicit but it is sufficient to estimate
$\sup_{z \in B_r} |\frac{\partial^{|\alpha|} (\xi)}{\partial
x^{\alpha}}(z)|$ like in (\ref{eqn:estimation-inverse}) for
$|\alpha|>1$ (note that in this case, we have
$\frac{\partial^{|\alpha|} (\Psi_i)}{\partial
x^{\alpha}}=\frac{\partial^{|\alpha|} (\xi_i)}{\partial
x^{\alpha}}$). Now we have to study the case $|\alpha|=1$. In this
case, writing the Jacobian matrix, we have
$$
1 + \frac{\partial \xi}{\partial x} = (1 + \frac{\partial
\chi}{\partial x})^{-1}\circ\Phi\,.
$$
Denoting by $|\!|\!|\, |\!|\!|$ the standard norm of linear
operators on a finite dimensional vector space we can assume that
$|\!|\!|\frac{\partial \chi}{\partial x}|\!|\!|<1$. Then, since
$(1 + \frac{\partial \chi}{\partial x})^{-1}=1+\sum_{q\geq
1}\big(\frac{\partial \chi}{\partial x}\big)^q$, we obtain
$$
\frac{\partial \xi}{\partial x} =\big( \sum_{q\geq
1}\big(\frac{\partial \chi}{\partial x}\big)^q\big)\circ\Phi\,.
$$
We then get
$$
|\!|\!|\frac{\partial \xi}{\partial x}|\!|\!| \leq M
|\!|\!|\frac{\partial \chi}{\partial x}|\!|\!|
$$
where $M$ is a positive constant and we conclude using the
equivalence of the norms. \QED

\section{Proof of Theorem \ref{thm:LeviPoissonS}}
\label{section:proof}

In order to prove Theorem \ref{thm:LeviPoissonS}, we will
construct by recurrence a sequence of local smooth coordinate
systems $(x^d,y^d) :=
(x_1^d,\hdots,x_{m}^d,y_1^d,\hdots,y_{n-m}^d)$, where
$(x^0,y^0)=(x_1,\hdots,x_{m},y_1,\hdots,y_{n-m})$ is the original
linear coordinate system as chosen in Section 2, which converges
to a local coordinate system
$(x^\infty,y^\infty)=(x_1^\infty,\hdots,x_{m}^\infty,y_1^\infty,\hdots,y_{n-m}^\infty)$,
in which the Poisson structure $\Pi$ has the required form.

For simplicity of exposition, we will assume that $\Pi$ is
$C^\infty$-smooth. However, in every step of the proof of Theorem
\ref{thm:LeviPoissonS}, we will only use differentiability of
$\Pi$ up to some finite order, and that's why our proof will also
work for finitely (sufficiently highly) differentiable Poisson
structures.

We will denote by $\Theta_d$ the local diffeomorphisms of
$(\bbR^n,0)$ such that
\begin{equation}
(x^d,y^d) (z) = (x^0,y^0) \circ \Theta_d (z)\,,
\end{equation}
where $z$ denotes a point of $(\bbR^n,0)$.

Denote by $\Pi^d$ the Poisson structure obtained from $\Pi$ by the
action of $\Theta_d$:
\begin{equation}
\Pi^d = (\Theta_d)_* \Pi\,.
\end{equation}

Of course, $\Pi^0 = \Pi$. Denote by $\{.,.\}_d$ the Poisson
bracket with respect to the Poisson structure $\Pi^d$. Then we
have
\begin{equation}
\{F_1,F_2\}_d(z) = \{F_1 \circ {\Theta_d},F_2 \circ {\Theta_d} \}
({\Theta_d}^{-1}(z))\,.
\end{equation}

Assume that we have constructed $ (x^d,y^d) = (x,y) \circ
\Theta_d$. Let us now construct $(x^{d+1},y^{d+1}) = (x,y) \circ
\Theta_{d+1}$. This construction consists of two steps : 1) find
an ``almost Levi factor'', i.e. coordinates $x^{d+1}_i$ such that
the error terms $\{x_i^{d+1},x_j^{d+1}\} - \sum_k c_{ij}^k
x_k^{d+1}$ are small, and 2) ``almost linearize'' it, i.e. find
the remaining coordinates $y^{d+1}$ such that in the coordinate
system $(x^{d+1}, y^{d+1})$ the Hamiltonian vector fields of the
functions $x^{d+1}_i$ are very close to linear ones. In fact, we
will define a local diffeomorphism $\theta_{d+1}$ of $(\R^n,0)$
and then put $\Theta_{d+1}=\theta_{d+1}\circ\Theta_d$. In
particular, we will have $\Pi^{d+1}=(\theta_{d+1})_* \Pi^d$ and
$(x^{d+1},y^{d+1})=(x^d,y^d)\circ(\Theta_d)^{-1}\circ\theta_{d+1}\circ\Theta_d$.

We write the current error terms (that we want to make smaller by
going from $(x^d,y^d)$ to $(x^{d+1},y^{d+1})$) as follows:
\begin{equation}
f_{ij}^d(x,y) = \{x_i,x_j\}_d  -  \sum_{k=1}^m c_{ij}^k x_k ,
\end{equation}
and
\begin{equation}
g_{i\alpha}^d(x,y) =  \{x_i,y_\alpha\}_d - \sum_{\beta=1}^{n-m}
a_{i\alpha}^\beta y_\beta .
\end{equation}

Consider the 2-cochain
\begin{equation}
f^d=\sum_{ij}f^d_{ij}\otimes\xi_i^*\wedge\xi_j^*
\end{equation}
of the Chevalley-Eilenberg complex associated to the $\fg$-module
$\mathcal{C}_r$, where $r = r_d$ depends on $d$ and is chosen as
follows:
\begin{equation}
r_d = 1 + {1 \over d+1}.
\end{equation}
In particular, $r_0 = 2$, $r_d/r_{d+1} \sim 1 + {1 \over d^2}$,
and $\lim_{d \to \infty} r_d = 1$ is positive. This choice of
radii $r_d$ means in particular that we will be able to arrange so
that the Poisson structure $\Pi^d = (\Theta_d)_* \Pi$ is defined
in the closed ball of radius $r_d$. (For this to hold, we will
have to assume that $\Pi$ is defined in the closed ball of radius
2, and show by recurrence that $B_{r_{d}} \subset \theta_d
(B_{r_{d-1}})$ for all $d \in \bbN$).

Put
\begin{equation}
\varphi^{d+1} := \sum_i \varphi_i^{d+1}\otimes\xi_i^* =
-S(t_d)\big( h(f^d)\big)\,,
\end{equation}
where $h$ is the homotopy operator as given in Lemma
\ref{lem:HomotopyOp}, $S$ is the smoothing operator and the
parameter $t_d$ is chosen as follows: take a real constant $t_0>1$
(which later on will be assumed to be large enough) and define the
sequence $(t_d)_{d\geq 0}$ by $t_{d+1}=t_d^{3/2}$. In other words,
we have
\begin{equation}
t_d=\exp\left(\left(\frac{3}{2}\right)^d\ln t_0\right)\;,\;\ln
t_0>0\,.
\end{equation}
The above choice of smoothing parameter $t_d$ is a standard one in
problems involving the Nash-Moser method, see, e.g.,
\cite{Hamilton-Complex1977,Hamilton-NashMoser1982}. The number $3
\over 2$ in the above formula is just a convenient choice. The
main point is that this number is greater than 1 (so we have a
very fast increasing sequence) and smaller than 2 (where 2
corresponds to the fact that we have a fast convergence algorithm
which ``quadratizes'' the error term at each step, i.e. go from an
``$\epsilon$-small'' error term to an ``$\epsilon^2$-small'' error
term).

According to Inequality (\ref{eqn:estimate-h1}), in order to
control the $C^k$-norm of $h(f^d)$ we need to control the
$C^{k+s}$-norm of $f^d$, i.e. we face a ``loss of
differentiability''. That's why in the above definition of
$\varphi^{d+1}$ we have to use the smoothing operator $S$, which
will allow us to compensate for this loss of differentiability.
This is a standard trick in the Nash-Moser method.

Next, consider the 1-cochains
\begin{eqnarray}
g^d &=& \sum_i\big(\sum_\alpha
g_{i\alpha}^d\frac{\partial}{\partial y_\alpha}\big)\otimes
\xi_i^* , \\
{\hat g}^d &=& g^d- \sum_i\big( \sum_\alpha
\{h(f^d)_i,y_\alpha\}_d \frac{\partial}{\partial
y_\alpha}\big)\otimes\xi_i^*
\end{eqnarray}
of the differential of the Chevalley-Eilenberg complex associated
to the $\fg$-module $\mathcal{Y}_r$, where $r = r_d = 1 + {1 \over
d+1}$, and put
\begin{equation}
\psi^{d+1} := \sum_\alpha
\psi^{d+1}_\alpha\frac{\partial}{\partial y_\alpha} = -S(t_d)\big(
h({\hat g}^d)\big)\,,
\end{equation}
where $h$ is the homotopy operator as given in Lemma
\ref{lem:HomotopyOp}, and $S(t_d)$ is the smoothing operator (with
the same $t_d$ as in the definition of $\varphi^{d+1}$).

Now define $\theta_{d+1}$ to be a local diffeomorphism of $\R^n$
given by
\begin{equation}
\theta_{d+1}:=Id+\chi^{d+1}:=Id+(\varphi^{d+1},\psi^{d+1}) \, ,
\end{equation}
where $(\varphi^{d+1},\psi^{d+1})$ now means
$(\varphi^{d+1}_1,\dots,\varphi^{d+1}_m,\psi^{d+1}_1,\dots,\psi^{d+1}_{n-m})$.
This finishes our construction of $\Theta_{d+1}=\theta_{d+1}\circ
\Theta_d$ and $(x^{d+1},y^{d+1})=(x,y)\circ\Theta_{d+1}$. This
construction is very similar to the analytic case
\cite{Zung-Levi2002}, except mainly for the use of the smoothing
operator. Another difference is that, for technical reasons, in
the smooth case considered in this paper we use the original
coordinate system and the transformed Poisson structures $\Pi^d$
for determining the error terms, while in the analytic case the
original Poisson structure and the transformed coordinate systems
are used. (In particular, the closed balls used in this paper are
always balls with respect to the original coordinate system --
this allows us to easily compare the Sobolev norms of functions on
them, i.e. bigger balls correspond to bigger norms).


In order to show that the sequence of diffeomorphisms defined
above converges to a smooth local diffeomorphism $\Theta_\infty$
and that the limit Poisson structure $({\Theta_\infty})_*\Pi$ is
in Levi normal form, we will have to control the norms of $\delta
f$ and ${\delta}{\hat g}^d$, where $\delta$ denotes the
differential of the corresponding Chevalley-Eilenberg complexes.
This will be done with the help of the following two simple
lemmas:

\begin{lem} \label{lem:df}
For every $i,j$ and $k$, we have
\begin{equation}
\delta f^d(\xi_i\wedge\xi_j\wedge\xi_k) = \oint_{ijk} \big( \sum_u
f^d_{iu}\frac{\partial f^d_{jk}}{\partial x_u}+ \sum_\alpha
g^d_{i\alpha}\frac{\partial f^d_{jk}}{\partial y_\alpha}\big)\,,
\end{equation}
where $\oint$ denotes the cyclic sum.
\end{lem}

\begin{lem} \label{lem:dg}
For every $i,j$ and $\alpha$, the coefficient of $\frac{\partial
}{\partial y_\alpha}$ in $\delta {\hat {g}}^d(\xi_i\wedge\xi_j)$
is
\begin{eqnarray*}
&-&\sum_u f_{iu}^d\frac{\partial g_{j\alpha}^d}{\partial x_u}-
\sum_\beta g_{i\beta}^d\frac{\partial g_{j\alpha}^d}{\partial
y_\beta}+ \sum_u f_{ju}^d\frac{\partial g_{i\alpha}^d}{\partial
x_u}+
\sum_\beta g_{j\beta}^d\frac{\partial g_{i\alpha}^d}{\partial y_\beta}\\
&+& \sum_u f_{iu}^d\frac{\partial
\{h(f^d)_j,y_\alpha\}_d}{\partial x_u} +\sum_\beta
g_{i\beta}^d\frac{\partial \{h(f^d)_j,y_\alpha\}_d}{\partial
y_\beta}\\
&-& \sum_u f_{ju}^d\frac{\partial
\{h(f^d)_i,y_\alpha\}_d}{\partial x_u}
-\sum_\beta g_{j\beta}^d\frac{\partial \{h(f^d)_i,y_\alpha\}_d}{\partial y_\beta}\\
&+& \{h(f^d)_i,g_{j\alpha}^d\}_d - \{h(f^d)_j,g_{i\alpha}^d\}_d\\
&+& \{y_\alpha,\sum_u f_{iu}\frac{\partial h(f^d)_j}{\partial
x_u}- \sum_u f_{ju}\frac{\partial h(f^d)_i}{\partial x_u} +
\sum_\beta g_{i\beta}\frac{\partial h(f^d)_j}{\partial y_\beta}-
\sum_\beta g_{j\beta}\frac{\partial h(f^d)_i}{\partial y_\beta}
\}\\
&-& \{y_\alpha,h(\delta f^d)_{ij}\} \ .
\end{eqnarray*}
\end{lem}

The first lemma is a direct consequence of the Jacobi identity
$\{x_i,\{x_j,x_k\}_d\}_d+\{x_j,\{x_k,x_i\}_d\}_d+\{x_k,\{x_i,x_j\}_d\}_d=0$.
The second one follows from the Jacobi identity
$\{x_i-h(f^d)_i,\{x_j-h(f^d)_j,y_\alpha\}_d\}_d+
\{x_j-h(f^d)_j,\{y_\alpha,x_i-h(f^d)_i\}_d\}_d+
\{y_\alpha,\{x_i-h(f^d)_i,x_j-h(f^d)_j\}_d\}_d=0$ and the homotopy
relation (\ref{eqn:homotopy'}). \QED

Roughly speaking, the above lemmas say that $\delta f^d$ and
${\delta}{\hat g}^d$ are ``quadratic functions'' in $f^d$, $g^d$
and their first derivatives, so if $f^d$ and $g^d$ are
``$\epsilon$-small'' then $\delta f^d$ and ${\delta}{\hat g}^d$
are ``$\epsilon^2$-small''.

Let us now give some expressions for the new error terms, which
will allow us to estimate their norms. Recall that the new error
terms after Step $d$ are
\begin{eqnarray}
f^{d+1}_{ij}(x,y) &=& \{x_i,x_j\}_{d+1}-\sum_k c_{ij}^k x_k \ , \\
g^{d+1}_{i\alpha}(x,y) &=& \{x_i,y_\alpha\}_{d+1}- \sum_\beta
a_{i\alpha}^\beta y_\beta \ .
\end{eqnarray}

We can also write, for instance,
\begin{equation}
f^{d+1}_{ij}(x,y)=
[\{x_i+\varphi^{d+1}_i,x_j+\varphi^{d+1}_j\}_d-\sum_k c_{ij}^k
(x_k+\varphi^{d+1}_k)](\theta_{d+1}^{-1}(x,y))\,.
\end{equation}

A simple direct computation shows that
\begin{eqnarray}
f^{d+1}_{ij} &=& \big[ (\delta \varphi^{d+1})_{ij}+
f^d_{ij}+Q^d_{ij}\big]\circ(\theta_{d+1})^{-1} \label{eqn:defif1} \ , \\
g^{d+1}_{i\alpha} &=& \big[ (\delta \psi^{d+1})_{i\alpha}+{\hat
{g}}^d_{i\alpha}+ T^d_{i\alpha}+ U^d_{i\alpha}
\big]\circ(\theta_{d+1})^{-1}\label{eqn:defig1} \ ,
\end{eqnarray}
where $Q^d_{ij}$ and $T^d_{i\alpha}$ are ``quadratic functions'',
namely
\begin{equation}
Q^d_{ij}=\sum_u(f^d_{iu}\frac{\partial \varphi_j^{d+1}}{\partial
x_u} -f^d_{ju}\frac{\partial \varphi_i^{d+1}}{\partial x_u}) +
\sum_\beta(g^d_{i\beta}\frac{\partial \varphi_j^{d+1}}{\partial
y_\beta} -g^d_{j\beta}\frac{\partial \varphi_i^{d+1}}{\partial
y_\beta}) + \{\varphi_i^{d+1},\varphi_j^{d+1}\}_d\,,
\label{eqn:defiQ}
\end{equation}
and
\begin{equation}
T^d_{i\alpha}=\sum_uf^d_{iu}\frac{\partial
\psi^{d+1}_\alpha}{\partial x_u} + \sum_\beta
g^d_{i\beta}\frac{\partial \psi^{d+1}_\alpha}{\partial y_\beta}
+\{\varphi^{d+1}_i,\psi_\alpha^{d+1}\}_d\,, \label{eqn:defiT}
\end{equation}
and $U^d_{i\alpha}$ is defined by
\begin{equation}
U^d_{i\alpha}=\{h(f^d)_i-S(t_d)h(f^d)_i,y_\alpha\}_d\,.
\label{eqn:defiU}
\end{equation}
Putting
\begin{eqnarray*}
Q^d &=& \sum_{ij} Q^d_{ij}\otimes\xi_i^*\wedge\xi_j^* \ ,\\
T^d &=& \sum_i(\sum_\alpha T^d_{i\alpha}\frac{\partial}{\partial
y_\alpha})\otimes\xi_i^* \ ,\\
U^d &=& \sum_i(\sum_\alpha U^d_{i\alpha}\frac{\partial}{\partial
y_\alpha})\otimes\xi_i^* \ ,
\end{eqnarray*}
we can write
\begin{eqnarray}
f^{d+1} &=& (\delta\varphi^{d+1}+f^d+Q^d)\circ (\theta_{d+1})^{-1}\ , \\
g^{d+1} &=& (\delta \psi^{d+1}+{\hat {g}}^d+T^d+U^d)\circ
(\theta_{d+1})^{-1}\,.
\end{eqnarray}

Equality (\ref{eqn:homotopy'}) allows us to give another
expression for $f^{d+1}$ and $g^{d+1}$, which will be more
convenient:

\begin{eqnarray}
f^{d+1} &=& \big[ \delta\big(\varphi^{d+1}+h(f^d)\big)+h(\delta
f^d)+Q^d\big]\circ
(\theta_{d+1})^{-1} \label{eqn:defif2}\ , \\
g^{d+1} &=& \big[ \delta \big(\psi^{d+1}+{h}({\hat {g}}^d)\big)+
{h}(\delta{\hat {g}}^d)+T^d+U^d\big]\circ
(\theta_{d+1})^{-1}\,.\label{eqn:defig2}
\end{eqnarray}

The following two technical lemmas about the norms will be the key
points of the proof of Theorem \ref{thm:LeviPoissonS}. In order to
formulate them, we need to introduce some positive constants $A,
l$ and $L$. Recall that we denote $s = [{n \over 2}] + 1$ (this
number $s$ appears in the Sobolev inequality and measures the
``loss of differentiability'' in our algorithm). Put $A = 6s + 9$.
We will use the fact that
\begin{equation}
A > 6s + 8 \ .
\end{equation}
Choose an auxilliary positive constant $\varepsilon < 1$ such that
\begin{equation}
-(1-\varepsilon)+A\varepsilon < -\frac{3}{4}\,. \label{eqn:ba}
\end{equation}
Choose an integer $l>s$ such that
\begin{eqnarray}
\frac{3s+5}{l-1} &<& \varepsilon \label{eqn:s}
\end{eqnarray}
(this is the number $l$ which appears in the formulation of
Theorem \ref{thm:LeviPoissonS}), and put
\begin{equation}
L=2l-1 \ .
\end{equation}
Recall also that $t_0>1$, $t_d=\exp ((3/2)^d\ln t_0)$ and
$r_d=1+\frac{1}{d+1}$ (note that we have
$r_{d+1}=r_d(1-\frac{1}{(d+2)^2})$). By choosing $t_0$ large
enough, we can assume that $t_d^{-1/2}<\frac{1}{(d+2)^2}$ for
every $d$.

\begin{lem}
\label{lem:Hamilton1} Suppose that $\Pi$ is defined on $B_{r_0}$
and satisfies the following  inequalities:
\begin{equation}
\| f^0\|_{l,r_0} < t_0^{-1}, \ \ \| g^0\|_{l,r_0} < t_0^{-1}, \ \
\|\Pi\|_{L,r_0}<t_0^A,\ \ \| f^0\|_{L,r_0} < t_0^{A}, \ \ \|
g^0\|_{L,r_0} < t_0^{A},
\end{equation}
where $t_0 > 1$ is a sufficiently large number. Then for every
nonnegative integer $d$, $\Pi^d$ is well-defined on $B_{r_d}$ and
we have the following estimates:

($1_d$)\ \ \  $\norm \chi^{d+1} \norm_{l,r_d} < t_d^{-1/2}$
(recall that $\chi^{d+1} =
-(\varphi^{d+1}_1,\hdots,\psi^{d+1}_{n-m})$)

($2_d$)\ \ \ $\norm\Pi^d\norm_{L,r_d}<t_d^A$

($3_d$)\ \ \  $\|\Pi^d\|_{l,r_d}<C {{d+1} \over {d+2}}$, where $C$
is a positive constant independent of $d$.

($4_d$)\ \ \ $\norm f^d\norm_{L,r_d} < t_d^{A}$ and $\norm
g^d\norm_{L,r_d} < t_d^{A}$

($5_d$)\ \ \ $\norm f^d\norm_{l,r_d} < t_d^{-1}$ and $\norm
g^d\norm_{l,r_d} < t_d^{-1}$
\end{lem}

Roughly speaking, Inequality ($1_d$) is the one which ensures the
convergence of $\Theta_d$ when $d \to \infty$ in $C^l$-topology.
Inequality ($3_d$) says that $\|\Pi^d\|_l$ stays bounded.
Inequality ($5_d$) means that the error terms converge to 0 very
fast in $C^l$-topology, while Inequalities ($2_d$) and ($4_d$)
mean that things don't ``get bad'' too fast in $C^L$-topology.

\begin{lem} Suppose that for an integer $k\geq l$, there exists a constant $C_k>0$ and
an integer $d_k\geq 0$ such that for any $d\geq d_k$, the
following inequalities are satisfied:
\begin{multline}
\| f^d\|_{k,r_d} < C_k t_d^{-1}, \ \ \| g^d\|_{k,r_d} < C_k
t_d^{-1}, \ \
\| f^d\|_{2k-1,r_d} < C_kt_d^{A}, \\
 \| g^d\|_{2k-1,r_d} < C_k t_d^{A},\ \ \|\Pi^d\|_{2k-1,r_d}< C_k t_d^{A},
\ \ \|\Pi^d\|_{k,r_d}< C_k(1- {1 \over d+2}).
\end{multline}
Then there exists a constant $C_{k+1}>0$ and an integer $d_{k+1} >
d_k$ such that, for any $d\geq d_{k+1}$, we have
\begin{itemize}
\item[i)]  $\norm \chi^{d+1}\norm_{k+1,r_d} < C_{k+1} t_d^{-1/2}$

\item[ii)]  $\norm\Pi^d\norm_{k+1,r_d}<C_{k+1}(1- {1 \over d+2})$

\item[iii)] $\norm f^d\norm_{k+1,r_d} < C_{k+1} t_d^{-1}$ and
$\norm g^d\norm_{k+1,r_d} < C_{k+1} t_d^{-1}$

\item[iv)] $\norm f^d\norm_{2k+1,r_d} < C_{k+1} t_d^A$,  $\norm
g^d\norm_{2k+1,r_d} < C_{k+1} t_d^A$ and
$\norm\Pi^d\norm_{2k+1,r_d}<C_{k+1}t_d^A$
\end{itemize}
\label{lem:Hamilton2}
\end{lem}

The above two technical lemmas will be proved in Section 5. Let us
now finish the proof of Theorem \ref{thm:LeviPoissonS} modulo
them.

{\bf Proof of Theorem \ref{thm:LeviPoissonS}}. Assume for the
moment that $\Pi$ is sufficiently close to its linear part, more
precisely, that the conditions of Lemma \ref{lem:Hamilton1} are
satisfied. Let $p$ be a natural number greater or equal to $l$
such that $\Pi$ is at least $C^{2p-1}$-smooth. Applying Lemma
\ref{lem:Hamilton1} to $\Pi$, and then applying Lemma
\ref{lem:Hamilton2} repetitively, we get the following inequality:
there exist an integer $d_p$ and a positive constant $C_p$ such
that for every $d\geq d_p$ we have
\begin{equation}
\|\chi^{d+1}\|_{p,r_d} \leq C_p t_d^{-1/2} = C_p \exp \Big(- {1
\over 2} \Big({3 \over 2}\Big)^d \ln t_0 \Big).
\end{equation}
The right hand side of the above inequality tends to 0
exponentially fast when $d \rightarrow \infty$. This, together
with Lemmas \ref{lem:composition} and \ref{lem:inverse}, implies
that
\begin{equation}
(\Theta_d)^{-1} = (\theta_1)^{-1} \circ \hdots \circ
(\theta_d)^{-1},
\end{equation}
where $\theta_d = {\rm Id} + \chi^d$, converges in $C^p$-topology
on the ball $B_{1}$ of radius 1 (we show in Lemma
\ref{lem:Hamilton1} that $(\Theta_d)^{-1}$ is well-defined on the
ball of radius $r_d > 1$). The fact that $\Theta_\infty = \lim_{d
\to \infty} \Theta_d$ is a local $C^p$-diffeomorphism should now
be obvious. It is also clear that $\Pi^\infty={(\Theta_\infty)}_*
\Pi$ is in Levi normal form. (Inequalities in $(5_d)$ of Lemma
\ref{lem:Hamilton1} measure how far is $\Pi^d$ from a Levi normal
form; these estimates tend to 0 when $d \to \infty$).

If $\Pi$ does not satisfy the conditions of Lemmas
\ref{lem:Hamilton1} and \ref{lem:Hamilton2}, then we may use the
following homothety trick: replace $\Pi$ by $\Pi^t = {1 \over t}
G(t)_* \Pi$ where $G(t): z \mapsto tz$ is a homothety, $t > 0$.
The limit $\lim_{t \to \infty} \Pi^t$ is equal to the linear part
of $\Pi$. So by choosing $t$ high enough, we may assume that
$\Pi^t$ satisfies the conditions of Lemmas \ref{lem:Hamilton1} and
\ref{lem:Hamilton2}. If $\Phi$ is a required local diffeomorphism
(coordinate transformation) for $\Pi^t$, then $G(1/t) \circ \Phi
\circ G(t)$ will be a required local smooth coordinate
transformation for $\Pi$. \QED

\section{Proof of the technical lemmas}
\label{section:technical}


{\it Proof of Lemma \ref{lem:Hamilton1}:} We prove this lemma by
induction on $d$. The main tools used are Leibniz-type
inequalities, and interpolation inequalities (Inequality
(\ref{eqn:interpolation})) involving $C^l$-norms, $C^L$-norms and
the norms in between. Roughly speaking, ($2_d$) and ($4_d$) will
follow from Leibniz-type inequalities. The proof of ($1_d$) and
($5_d$) will make substantial use of interpolation inequalities.
Point ($3_d$) will from from an analog of ($1_d$) and Leibniz-type
inequalities.


In order to simplify the notations, we will use the letter $M$ to
denote a constant, {\it which does not depend on $d$} but which
varies from inequality to inequality (i.e. it depends on the line
where it appears).

We begin the reduction at $d=0$. For $d=0$, the only point to be
checked is $(1_0)$. We will use a property of the smoothing
operator (equation (\ref{eqn:smoothing1})), the estimate of the
homotopy operator (\ref{eqn:estimate-h1}) and the interpolation
relation (\ref{eqn:interpolation}).

Recall that $\varphi^1=-S(t_0)\big(h(f^0)\big)$. We then have
\begin{eqnarray*}
\|\varphi^1\|_{l,r_0} &\leq& M \|h(f^0)\|_{l,r_0} \quad {\mbox { by }} (\ref{eqn:smoothing1})  \\
                      &\leq& M \|f^0\|_{l+s,r_0} \quad {\mbox { by }} (\ref{eqn:estimate-h1})\\
                      &\leq& M \|f^0\|_{l,r_0}^{\frac{l-s-1}{l-1}}
                               \|f^0\|_{L,r_0}^{\frac{s}{l-1}}
                               \quad {\mbox { by }} (\ref{eqn:interpolation})\\
                      &\leq& M t_0^{-\frac{l-s-1}{l-1}}
                               t_0^{A\frac{s}{l-1}}
\end{eqnarray*}
On the other hand, we have $\psi^1=-S(t_0)\big(h({\hat
g}^0)\big)$, then
\begin{eqnarray*}
\|\psi^1\|_{l,r_0} &\leq& M\|{\hat g}^0\|_{l+s,r_0} \quad {\mbox {
by }} (\ref{eqn:smoothing1})
{\mbox { and }} (\ref{eqn:estimate-h1}) \\
                 &\leq& M\|g^0+\{h(f^0),y\}_0\|_{l+s,r_0}\\
                 &\leq& M(\|g^0\|_{l+s,r_0}+\|\Pi^0\|_{l+s,r_0}\|h(f^0)\|_{l+s+1,r_0})\,;
\end{eqnarray*}
note that since the definition of ${\hat g^0}$ involves the first
derivatives of $h(f^0)$ we have to estimate by
$\|h(f^0)\|_{l+s+1}$. Now, using (\ref{eqn:estimate-h1}) and the
interpolation relation (\ref{eqn:interpolation}), we get
\begin{eqnarray*}
\|\psi^1\|_{l,r_0} &\leq& M(\|g^0\|_{l+s,r_0}+\|\Pi^0\|_{l+s,r_0}\|f^0\|_{l+2s+1,r_0})\\
                   &\leq& M(\|g^0\|_{l,r_0}^{\frac{l-s-1}{l-1}}\|g^0\|_{L,r_0}^{\frac{s}{l-1}}
                 +\|\Pi^0\|_{l,r_0}^{\frac{l-s-1}{l-1}}\|\Pi^0\|_{L,r_0}^{\frac{s}{l-1}}
                 \|f^0\|_{l,r_0}^{\frac{l-2s-2}{l-1}}\|f^0\|_{L,r_0}^{\frac{2s+1}{l-1}})\\
                 &\leq& M(t_0^{-\frac{l-s-1}{l-1}+A\frac{s}{l-1}}
                     +t_0^{-\frac{l-2s-2}{l-1}+A\frac{2s+1}{l-1}+A\frac{s}{l-1}})\,.
\end{eqnarray*}
Since $\frac{s}{l-1}<\frac{2s+1}{l-1}$ we have
$-\frac{l-s-1}{l-1}+A\frac{s}{l-1}<-\frac{l-2s-2}{l-1}+A\frac{2s+1}{l-1}+A\frac{s}{l-1}$.
Therefore, we obtain
$$
\|\chi^1\|_{l,r_0}<M
t_0^{-\frac{l-2s-2}{l-1}+A\frac{3s+1}{l-1}}\,.
$$
By assumptions (\ref{eqn:s}), we see that $\frac{2s+1}{l-1}$ and
$\frac{3s+1}{l-1}$ are strictly smaller than $\varepsilon$.
Therefore, $-\frac{l-2s-2}{l-1}+A\frac{3s+1}{l-1}$ is strictly
smaller than $-(1-\varepsilon)+A\varepsilon$. Then, according to
inequality (\ref{eqn:ba}), we have $\|\chi^1\|_{l,r_0}\leq
Mt_0^{-\mu}$ with $-\mu<-3/4<-1/2$. We may choose $t_0$
sufficiently large such that $Mt_0^{-\mu}<t_0^{-1/2}$, which gives
\begin{equation}
\|\chi^1\|_{l,r_0} < t_0^{-1/2}\,. \label{eqn:chi_l}
\end{equation}

 Now, by induction, we suppose that for some $d\geq 0$, $\Pi^d$ is well
 defined on $B_{r_d}$ and that the
inequalities $(1_d),\hdots,(5_d)$ are true. We will show that they
still hold when we replace $d$ by $d+1$. To simplify the writing
we will omit the index $r_d$ in the norms, unless the radius in
question is different from $r_d$.

 Since $r_{d+1}=r_d(1-\frac{1}{(d+2)^2})$, according to Inequality
$(1_d)$ ($\|\chi^{d+1}\|_{l,r_d}<t_d^{-1/2}$ which is, by
assumption, strictly smaller than $\frac{1}{(d+2)^2}$ for every
$d$) and Lemma \ref{lem:boules}, we know that $B_{r_{d+1}}$ is
included in $\theta_{d+1}(B_{r_d})$ and so $\Pi^{d+1}$ will be
well defined on $B_{r_{d+1}}$.

$\bullet$ Proof of $(1_{d+1})$: absolutely similar to the proof of
($1_0$) given above.

$\bullet$ Proof of $(2_{d+1})$: Recall that, due to the fact that
$\Pi^{d+1} = (\Theta_{d+1})_*\Pi = (\theta_{d+1})_*\Pi^d$, we have
$$
\{x_i,x_j\}_{d+1} = \{x_i + \varphi^{d+1}_i,x_j +
\varphi^{d+1}_j\}_d \circ (\theta_{d+1})^{-1},
$$
and similar formulas for  ${\{ x_i,y_\alpha \}}_{d+1}$ and
${\{y_\alpha,y_\beta\}}_{d+1}$.


Applying Lemmas \ref{lem:compositionL} and \ref{lem:inverse} we
obtain
\begin{eqnarray*}
\|{\{ x_i,x_j\}}_{d+1} \|_{L,r_{d+1}} &\leq&
\|\{x_i+\varphi_i^{d+1},x_j+\varphi^{d+1}_j\}_d\|_{L,r_d}
P(\|\chi^{d+1}\|_{l,r_d})\\
    &+&
    \|\{x_i+\varphi_i^{d+1},x_j+\varphi^{d+1}_j\}_d\|_{l,r_d}
    \|\chi^{d+1}\|_{L,r_d}Q(\|\chi^{d+1}\|_{l,r_d})\,,
\end{eqnarray*}
where $P$ and $Q$ are polynomials functions which do not depend on
$d$. By the Leibniz rule of derivation, the term
$\|\{x_i+\varphi_i^{d+1},x_j+\varphi^{d+1}_j\}_d\|_l$ may be
estimated by $M\|\Pi^d\|_l(1+\|\varphi^{d+1}\|_{l+1})^2$ and,
using the same technic as in the proof of $(1_0)$, we can write
$\|\varphi^{d+1}\|_{l+1}<t_d^{-1/2}$. Therefore, using $(3_d)$, we
can write $\|\{x_i+\varphi_i^{d+1},x_j+\varphi^{d+1}_j\}_d\|_l
\leq M$. Consequently, we have
$$
\|{\{ x_i,x_j\}}_{d+1} \|_{L,r_{d+1}} \leq M
\|\{x_i+\varphi_i^{d+1},x_j+\varphi^{d+1}_j\}_d\|_{L,r_d} +M
\|\chi^{d+1}\|_{L,r_d}\,.
$$
We first study the term $\chi^{d+1}$. Actually, we will estimate
$\|\chi^{d+1}\|_{L+1}$ rather than $\|\chi^{d+1}\|_{L}$ because it
will be useful for the estimation of
$\|\{x_i+\varphi_i^{d+1},x_j+\varphi^{d+1}_j\}_d\|_{L}$. We first
write $\|\varphi^{d+1}\|_{L+1}\leq M t_d^{s+1}\|h(f^d)\|_{L-s}$ by
the property (\ref{eqn:smoothing1}) of the smoothing operator.
Using the estimate (\ref{eqn:estimate-h1}) for the homotopy
operator, we obtain $\|\varphi^{d+1}\|_{L+1}\leq M
t_d^{s+1}\|f^d\|_{L} \leq M t_d^{A+s+1}$. Now, we have
\begin{eqnarray*}
\|\psi^{d+1}\|_{L+1} &\leq& M t_d^{3s+2}\|{\hat h}({\hat
g}^d)\|_{L-3s-1}\,
{\mbox{  by }}\, (\ref{eqn:smoothing1}) \\
                    &\leq& M t_d^{3s+2}\|{\hat g}^d\|_{L-2s-1} {\mbox{  by }}
\, (\ref{eqn:estimate-h1})
\end{eqnarray*}

Then the definition of ${\hat g}^d$, the Leibniz rule of
derivation (recall that $L=2l-1$) and Inequality
(\ref{eqn:estimate-h1}) give
\begin{eqnarray*}
\|\psi^{d+1}\|_{L+1} &\leq& M t_d^{3s+2}(\|g^d\|_{L-2s-1}
+\|\Pi^d\|_{L-2s-1}\|h(f^d)\|_{l-s-1+1}\\
                     & & + \|\Pi^d\|_{l-s-1}\|h(f^d)\|_{L-2s-1+1})\\
                    &\leq& M t_d^{3s+2}(\|g^d\|_{L}
+\|\Pi^d\|_{L}\|f^d\|_{l}+\|\Pi^d\|_{l}\|f^d\|_{L})\\
                    &\leq& M t_d^{A+3s+2}\,.
\end{eqnarray*}
Therefore, we can write
$$
\|\chi^{d+1}\|_{L+1,r_d} \leq M t_d^{A+3s+2}\,.
$$

Note that in the same way as in the proof of $(1_0)$, one can show
that $\|\chi^{d+1}\|_{l+1,r_d}< t_d^{-1/2}$ and then, using once
more the Leibniz formula of the derivation of a product, we get
\begin{eqnarray*}
\|\{x_i+\varphi_i^{d+1},x_j+\varphi^{d+1}_j\}_d\|_{L,r_d} &\leq&
M\big(\|\Pi^d\|_L
(1+\|\varphi^{d+1}\|_{l+1})^2\\
& & \quad +\|\Pi^d\|_l(1+\|\varphi^{d+1}\|_{l+1})(1+\|\varphi^{d+1}\|_{L+1})\big)\\
 &\leq& M(\|\Pi^d\|_L+\|\varphi^{d+1}\|_{L+1}+1)\\
 &\leq& M(\|\Pi^d\|_L+t_d^{A+3s+2}+1)\\
 &\leq& Mt_d^{A+3s+2}\,.
\end{eqnarray*}
Exactly in the same way, we can estimate the terms
$\|\{x_i+\varphi_i^{d+1},y_\alpha+\psi^{d+1}_\alpha\}_d\|_{L,r_d}$
and
$\|\{y_\alpha+\psi_\alpha^{d+1},y_\beta+\psi^{d+1}_\beta\}_d\|_{L,r_d}$
by $Mt_d^{A+3s+2}$. To conclude, since by our choice $A=6s+9$ we
have $A+3s+2<3A/2$, these estimates lead to
$\|\Pi^{d+1}\|_{L,r_{d+1}}\leq Mt_d^D$ where $D$ is a positive
constant such that $D<3A/2$. Therefore, we may choose $t_0$ large
enough (in a way which does not depend on $d$) in order to obtain
$\|\Pi^{d+1}\|_{L,r_{d+1}} < t_d^{3A/2}=t_{d+1}^A$.\\


$\bullet$ Proof of $(3_{d+1})$: Recall again that we have
$$
\{x_i,x_j\}_{d+1} = \{x_i + \varphi^{d+1}_i,x_j +
\varphi^{d+1}_j\}_d \circ (\theta_{d+1})^{-1},
$$
and similar formulas involving also $y_i$-components.

The estimates in Lemmas \ref{lem:composition} and
\ref{lem:inverse} give
\begin{equation}
\|\Pi^{d+1}\|_{l,r_{d+1}}\leq\|\Lambda^{d+1}\|_{l,r_d}(1+
P(\|\chi^{d+1}\|_{l,r_d})),
\end{equation}
 where $p$ is a polynomial (which does not depend on $d$) with vanishing
 constant term, and
\begin{multline}
\Lambda^{d+1} = \sum_{ij} \{x_i + \varphi^{d+1}_i,x_j +
\varphi^{d+1}_j\}_d {\partial \over
\partial x_i} \wedge {\partial\over \partial x_j} + \\ + \sum_{i\alpha}
\{x_i + \varphi^{d+1}_i,y_\alpha +\psi^{d+1}_\alpha\}_d {\partial
\over
\partial x_i} \wedge {\partial\over \partial y_\alpha}
+ \sum_{\alpha \beta}\{y_\alpha + \psi^{d+1}_\alpha,y_\beta +
\psi^{d+1}_\beta\}_d {\partial \over
\partial y_\alpha} \wedge {\partial\over \partial y_\beta}.
\end{multline}

Notice that $\Lambda^{d+1}$ is equal to $\Pi^d$ plus terms which
involve $\chi^{d+1}$ and the $\Pi^{d}$-bracket. Hence, by the
Leibniz formula, we can write
\begin{equation}
\|\Lambda^{d+1}\|_{l,r_d} \leq \|\Pi^{d}\|_{l,r_d}(1 + M
\|\chi^{d+1}\|_{l+1,r_d})^2,
\end{equation}
which implies that
\begin{equation}
\|\Pi^{d+1}\|_{l,r_{d+1}}\leq \|\Pi^{d}\|_{l,r_d}(1 + M
\|\chi^{d+1}\|_{l+1,r_d})^2 (1+p(\|\chi^{d+1}\|_{l,r_d})) .
\end{equation}

Similarly to the proof of $(1_{d+1})$, it is easy to see that
$\|\chi^{d+1}\|_{l+1,r_d}<t_d^{-1/2}$, which is exponentially
small when $d \rightarrow \infty$. By choosing the constant $t_0$
large enough, we may assume that
\begin{equation}
(1 + M \|\chi^{d+1}\|_{l+1,r_d})^2 (1+p(\|\chi^{d+1}\|_{l,r_d})) <
1 + {1 \over (d+1)(d+3)}
\end{equation}
Together with the induction hypothesis $\|\Pi^{d}\|_{l,r_d} <{
C(d+1) \over (d+2)}$, we get
\begin{equation}
\|\Pi^{d+1}\|_{l,r_{d+1}} < { C(d+1) \over (d+2)}(1 + {1 \over
(d+1)(d+3)}) = {C(d+2) \over (d+3)}.
\end{equation}

$\bullet$ Proof of $(4_{d+1})$: Recall that
$$
f^{d+1}_{ij}=\{x_i,x_j\}_{d+1}-\sum_k c_{ij}^k x_k\ .,
$$
It is easy to check that for every $i$ and $j$,
$$
\|\sum_k c_{ij}^k x_k\|_{L,r_{d+1}}\leq B \|\Pi\|_{1,r_{d+1}}\leq
B \|\Pi\|_{L,r_{d+1}}\,,
$$
where $B$ is a positive constant which only depends on the
dimension $n$. This implies immediately that
$$
\|f^{d+1}\|_{L,r_{d+1}} \leq (B+1) \|\Pi^{d+1}\|_{L,r_{d+1}} \,.
$$
In Point $(2_{d+1})$, we showed that
$\|\Pi^{d+1}\|_{L,r_{d+1}}<Mt_d^D$ where $D$ is a positive
constant such that $D<3A/2$ therefore, replacing $t_0$ by a larger
real number (which of course does not depend on $d$) if necessary,
we have $\|f^{d+1}\|_{L,r_{d+1}}<t_d^{3A/2}=t_{d+1}^A$. The
estimate of
$\|g^{d+1}\|_{L,r_{d+1}}$ can be done in the same way.\\

$\bullet$ Proof of $(5_{d+1})$ :

Recall the formula (\ref{eqn:defif2})
$$
f^{d+1}=\big[ \delta\big(\varphi^{d+1}+h(f^d)\big)+h(\delta
f^d)+Q^d\big]\circ (\theta_{d+1})^{-1}\,.
$$
We then have, using lemmas \ref{lem:composition} and
\ref{lem:inverse}
\begin{equation}
\|f^{d+1}\|_{l,r_{d+1}}\leq
M\|\delta\big(\varphi^{d+1}+h(f^d)\big) +h(\delta
f^d)+Q^d\|_{l,r_d}(1+ P(t_d^{-1/2}))\,, \label{eqn:compofl}
\end{equation}
where $P$ is a polynomial function.

Thus, we only have to estimate
$\|\delta\big(\varphi^{d+1}+h(f^d)\big)+h(\delta
f^d)+Q^d\|_{l,r_d}$. To do that, we use the second property of the
smoothing operator (\ref{eqn:smoothing2}), the estimate of the
homotopy operator (\ref{eqn:estimate-h1}) and the interpolation
inequality.

 We first write
\begin{eqnarray*}
\|\delta\big(\varphi^{d+1}+h(f^d)\big)\|_l &\leq& M\|h(f^d)-S(t_d)h(f^d)\|_{l+1}\\
 &\leq& M t_d^{-1}\|h(f^d)\|_{l+2}\quad {\mbox { by (\ref{eqn:smoothing2})}}\\
 &\leq& M t_d^{-1}\|f^d\|_{l+2+s}\quad {\mbox { by (\ref{eqn:estimate-h1})}}\\
 &\leq& M t_d^{-1}\|f^d\|_l^{\frac{l-s-3}{l-1}}\|f^d\|_L^{\frac{s+2}{l-1}}\quad
{\mbox {by (\ref{eqn:interpolation})}}
\end{eqnarray*}
Then, we have
\begin{equation}
\|\delta\big(\varphi^{d+1}+h(f^d)\big)\|_l \leq M
t_d^{-1-\frac{l-s-3}{l-1}+ A\frac{s+2}{l-1}}\,. \label{eqn:Q1l}
\end{equation}
Next, we write
\begin{eqnarray*}
\|h(\delta f^d)\|_l &\leq& M \|\delta f^d\|_{l+s} \quad {\mbox { by (\ref{eqn:smoothing1})}} \\
 &\leq& M (\|f^d\|_{l+s}\|f^d\|_{l+s+1}+\|g^d\|_{l+s}\|f^d\|_{l+s+1})\quad
{\mbox {by Lemma \ref{lem:df}}}\\
 &\leq& M (\|f^d\|_{l+s+1}^2+\|g^d\|_{l+s+1}\|f^d\|_{l+s+1})\\
 &\leq& M (\|f^d\|_l^{2\frac{l-s-2}{l-1}}\|f^d\|_L^{2\frac{s+1}{l-1}}+
\|g^d\|_l^{\frac{l-s-2}{l-1}}\|g^d\|_L^{\frac{s+1}{l-1}}\|f^d\|_l^{\frac{l-s-2}{l-1}}
\|f^d\|_L^{\frac{s+1}{l-1}})\,.
\end{eqnarray*}
Thus,
\begin{equation}
\|h(\delta f^d)\|_l \leq M
t_d^{-2\frac{l-s-2}{l-1}+2A\frac{s+1}{l-1}}\,. \label{eqn:Q2l}
\end{equation}
Finally, by the definition (\ref{eqn:defiQ}) of $Q^d$ , we have
$$\|Q^d\|_l\leq M (\|f^d\|_l\|\varphi^{d+1}\|_{l+1}+\|g^d\|_l\|\varphi^{d+1}\|_{l+1}
+\|\Pi^d\|_l\|\varphi^{d+1}\|_{l+1}^2)\,.$$ In the same way as in
the proof of the point $(1_0)$, we can easily show that
$\|\chi^{d+1}\|_{l+1}<M
t_d^{-\frac{l-2s-3}{l-1}+A\frac{3s+3}{l-1}}$. Therefore, we can
write
\begin{equation}
\|Q^d\|_l \leq M (t_d^{-1-\frac{l-2s-3}{l-1}+A\frac{3s+3}{l-1}}
+t_d^{-2\frac{l-2s-3}{l-1}+2A\frac{3s+3}{l-1}})\,. \label{eqn:Q3l}
\end{equation}
Combining (\ref{eqn:Q1l}), (\ref{eqn:Q2l}) and (\ref{eqn:Q3l}) we
obtain
$$
\|f^{d+1}\|_{l,r_{d+1}} < M
t_d^{-2\frac{l-2s-3}{l-1}+2A\frac{3s+3}{l-1}}\,.
$$

Now, by (\ref{eqn:s}),  $\frac{2s+2}{l-1}$ and $\frac{3s+3}{l-1}$
are strictly smaller than $\varepsilon$, and then
$-2\frac{l-2s-3}{l-1}+2A\frac{3s+3}{l-1}$ is strictly smaller than
$-2(1-\varepsilon)+2A\varepsilon$. To finish, the inequality
(\ref{eqn:ba}) gives $\|f^{d+1}\|_{l,r_{d+1}}<M t_d^{-\alpha}$
where $-\alpha<-\frac{3}{2}$. We may choose $t_0$ large enough (in
a way which depends on $\alpha$ but not on $d$) in order to obtain
$\|f^{d+1}\|_{l,r_{d+1}}<t_d^{-\frac{3}{2}}=t_{d+1}^{-1}$.\\


Now, we apply the same technic to estimate
$\|g^{d+1}\|_{l,r_{d+1}}$. Recall the formula (\ref{eqn:defig2})
$$
g^{d+1} = \big[ \delta\big(\psi^{d+1}+{h}({\hat {g}}^d)\big)+ {h}(
\delta{\hat {g}}^d)+T^d +U^d \big]\circ (\theta_{d+1})^{-1}\,.
$$
In the same way as above, according to Lemmas
\ref{lem:composition} and \ref{lem:inverse}, we just have to
estimate $\| \delta\big(\psi^{d+1}+{h}({\hat {g}}^d)\big)+
{h}({\delta}{\hat {g}}^d)+T^d +U^d\|_l$. We first write
\begin{eqnarray*}
\| \delta \big(\psi^{d+1}+{h}({\hat {g}}^d)\big)\|_l &\leq&
M \|-S(t_d){h}({\hat {g}}^d) +{h}({\hat {g}}^d)\|_{l+1}\\
 &\leq& M t_d^{-1}\|{h}({\hat {g}}^d)\|_{l+2} \quad {\mbox { by (\ref{eqn:smoothing2})}}\\
 &\leq& M t_d^{-1}\|{\hat {g}}^d\|_{l+s+2}\quad {\mbox { by (\ref{eqn:estimate-h1})}}\\
 &\leq& M t_d^{-1} \|g^d+{\{ h(f^d),y \}}_d \|_{l+s+2}\\
 &\leq& M t_d^{-1} (\|g^d\|_{l+s+2}+\|\Pi^d\|_{l+s+2}\|h(f^d)\|_{l+s+3})\\
 &\leq& M t_d^{-1} (\|g^d\|_{l+s+2}+\|\Pi^d\|_{l+s+2}\|f^d\|_{l+2s+3})
\end{eqnarray*}
Using the interpolation inequality (\ref{eqn:interpolation}), we
obtain
\begin{eqnarray*}
\| \delta \big(\psi^{d+1}+{h}({\hat {g}}^d)\big)\|_l &\leq&
 M t_d^{-1}
(\|g^d\|_l^{\frac{l-s-3}{l-1}}\|g^d\|_L^{\frac{s+2}{l-1}}\\
 & & \quad\quad +\|\Pi^d\|_l^{\frac{l-s-3}{l-1}}
\|\Pi^d\|_L^{\frac{s+2}{l-1}}\|f^d\|_l^{\frac{l-2s-4}{l-1}}\|f^d\|_L^{\frac{2s+3}{l-1}})\\
 &\leq & M (t_d^{-1-\frac{l-s-3}{l-1}+A\frac{s+2}{l-1}}+
t_d^{-1-\frac{l-2s-4}{l-1}+A\frac{3s+5}{l-1}})
\end{eqnarray*}

and then, since $\frac{s+2}{l-1}<\frac{2s+3}{l-1}$,
\begin{equation}
\|\delta \big(\psi^{d+1}+ h({\hat {g}}^d)\big)\|_l\leq M
t_d^{-1-\frac{l-2s-4}{l-1}+A\frac{3s+5}{l-1}}\,. \label{eqn:T1l}
\end{equation}

We also have, by the estimate of the homotopy operator
(\ref{eqn:estimate-h1}),
$$\| h(\delta {\hat {g}}^d)\|_l\leq M\| \delta {\hat {g}}^d\|_{l+s}\,,$$
and using Lemma \ref{lem:dg} and the interpolation inequality
(\ref{eqn:interpolation}), we obtain
\begin{eqnarray*}
\| h(\delta {\hat {g}}^d)\|_l & \leq &  M\Big( \|f^d\|_{l+s}
    \|g^d\|_{l+s+1}+\|g^d\|_{l+s}\|g^d\|_{l+s+1}  +
    \|\Pi^d\|_{l+s+1}\|f^d\|_{l+s}\|h(f^d)\|_{l+s+2} \\
& & +\|\Pi^d\|_{l+s+1}\|g^d\|_{l+s}\|h(f^d)\|_{l+s+2}
    +\|\Pi^d\|_{l+s}\|h(f^d)\|_{l+s+1}\|g^d\|_{l+s+1}\\
& & + \|\Pi^d\|_{l+s}\|f^d\|_{l+s+1}\|h(f^d)\|_{l+s+2}
    + \|\Pi^d\|_{l+s}\|g^d\|_{l+s+1}\|h(f^d)\|_{l+s+2}\\
& & + \|\Pi^d\|_{l+s}\|h(\delta f^d)\|_{l+s+1}\Big)\\
& \leq & M\Big( \|f^d\|_{l+2s+2}\|g^d\|_{l+2s+2}
       +\|g^d\|_{l+2s+2}^2 + \|\Pi^d\|_{l+s+1}\|f^d\|_{l+2s+2}^2 \\
& &  \quad + \|\Pi^d\|_{l+s+1}\|g^d\|_{l+2s+2}\|f^d\|_{l+2s+2}
           +\|\Pi^d\|_{l+s}\|h(\delta f^d)\|_{l+s+1}\Big) \\
&\leq& M \Big(t_d^{-2\frac{l-2s-3}{l-1} +
2A\frac{2s+2}{l-1}+A\frac{s+1}{l-1}} +\|\Pi^d\|_{l+s}\|h(\delta
f^d)\|_{l+s+1}\Big)
\end{eqnarray*}

In the same way as in the proof of (\ref{eqn:Q2l}) one can show
that
$$
\| h(\delta f^d)\|_{l+s+1}\leq M t_d^{-2\frac{l-2s-3}{l-1} +
2A\frac{2s+2}{l-1}}
$$

This gives, applying the interpolation inequality to
$\|\Pi^d\|_{l+s}$,
\begin{equation}
\| h( \delta {\hat {g}}^d)\|_l \leq M t_d^{-2\frac{l-2s-3}{l-1} +
2A\frac{2s+2}{l-1}+A\frac{s+1}{l-1}} \,. \label{eqn:T2l}
\end{equation}

Now, recalling the definition of $T$ (\ref{eqn:defiT}), we have
$$
\|T^d\|_l \leq
M(\|f^d\|_l+\|g^d\|_l+\|\Pi^d\|_l\|\varphi^{d+1}\|_{l+1})
\|\psi^{d+1}\|_{l+1}\,,
$$
and using the estimate of $\|\chi^{d+1}\|_{l+1,r_d}$ given above,
we can show that
\begin{equation}
\|T^d\|_l \leq M (t_d^{-1-\frac{l-2s-3}{l-1}+A\frac{3s+3}{l-1}}
+t_d^{-2\frac{l-2s-3}{l-1}+2A\frac{3s+3}{l-1}})\,. \label{eqn:T3l}
\end{equation}
Finally, by the definition of $U^d$ (\ref{eqn:defiU}), we can
write
\begin{eqnarray*}
\|U^d\|_{l} &\leq& M \|\Pi^d\|_l \|h(f^d)-S(t_d)h(f^d)\|_{l+1}\\
            &\leq& M \|\Pi^d\|_l t_d^{-1}\|h(f^d)\|_{l+2}
            \quad {\mbox { by (\ref{eqn:smoothing2})}}\\
            &\leq& M t_d^{-1} \|f^d\|_{l+s+2}
            \quad {\mbox { by ($3_d$) and (\ref{eqn:estimate-h1})}}\\
            &\leq& M t_d^{-1} \|f^d\|_l^{\frac{l-s-3}{l-1}}
            \|f^d\|_L^{\frac{s+2}{l-1}}
            \quad {\mbox { by (\ref{eqn:interpolation})}}\,.
\end{eqnarray*}
We then obtain
\begin{equation}
\|U^d\|_{l}\leq M t_d^{-1-\frac{l-s-3}{l-1}+A\frac{s+2}{l-1}}\,.
\label{eqn:T4l}
\end{equation}

Combining (\ref{eqn:T1l}), (\ref{eqn:T2l}), (\ref{eqn:T3l}) and
(\ref{eqn:T4l}), we obtain
$$
\|g^{d+1}\|_{l,r_{d+1}} < M
t_d^{-2\frac{l-2s-4}{l-1}+2A\frac{3s+5}{l-1}}\,,
$$
and we can conclude in the same way as for the estimate
of $\|f^{d+1}\|_{l,r_{d+1}}$.\\

Lemma \ref{lem:Hamilton1} is proved.\QED


{\it Proof of Lemma \ref{lem:Hamilton2}:} The main tools used in
the proof of this lemma are the same as in the previous lemma:
Leibniz-type inequalities and interpolation inequalities. To
simplify the notations, we will denote by $M_k$ a positive
constant which depends on $k$ {\it but not on $d$} and which
varies from inequality to inequality.

$\bullet$ Proof of {\it (i)}: If $d \geq d_k$, we have
\begin{eqnarray*}
\|\varphi^{d+1}\|_{k+1,r_d}= \|S(t_d)h(f^d)\|_{k+1} &\leq& M_k
\|f^d\|_{k+s+1}
\quad {\mbox { by }} (\ref{eqn:smoothing1}) {\mbox { and }} (\ref{eqn:estimate-h1}) \\
      &\leq& M_k \|f^d\|_k^{\frac{k-s-2}{k-1}} \|f^d\|_{2k-1}^{\frac{s+1}{k-1}}
\quad {\mbox { by }} (\ref{eqn:interpolation}) \\
      &\leq& M_k t_d^{-\frac{k-s-2}{k-1}+A\frac{s+1}{k-1}} \,.
\end{eqnarray*}
In the same way, we get
\begin{eqnarray*}
\|\psi^{d+1}\|_{k+1,r_d} = \|S(t_d){\hat h}({\hat g}^d)\|_{k+1}
      &\leq& M_k \|{\hat g}^d\|_{k+s+1}\\
      &\leq& M_k \|g^d+{\{h(f^d),y\}}_d\|_{k+s+1}\\
      &\leq& M_k (\|g^d\|_{k+s+1}+\|\Pi^d\|_{k+s+1}\|h(f^d)\|_{k+s+2})\\
      &\leq& M_k (\|g^d\|_{k+s+1}+\|\Pi^d\|_{k+s+1}\|f^d\|_{k+2s+2})\\
      &\leq& M_k (t_d^{-\frac{k-s-2}{k-1}+A\frac{s+1}{k-1}}+
                t_d^{A\frac{s+1}{k-1}-\frac{k-2s-3}{k-1}+A\frac{2s+2}{k-1}})
\end{eqnarray*}
Therefore, we have
$$
\|\chi^{d+1}\|_{k+1,r_d} \leq M_k
t_d^{-\frac{k-2s-3}{k-1}+A\frac{3s+3}{k-1}}\,.
$$
According to the inequality (\ref{eqn:s}), the terms
$\frac{2s+2}{k-1}$ and $\frac{3s+3}{k-1}$ are strictly smaller
than $\varepsilon$. Then, $-\frac{k-2s-3}{k-1}+A\frac{3s+3}{k-1}$
is strictly smaller than $-(1-\varepsilon)+A\varepsilon$.
Therefore, by (\ref{eqn:ba}), we can write
$\|\chi^{d+1}\|_{k+1,r_d}<M t_d^{-\mu}$ with $-\mu<-3/4<-1/2$. We
conclude that there exists a positive integer $d_{k+1} > d_k$ such
that $\forall d \geq d_{k+1}$,
$\|\chi^{d+1}\|_{k+1,r_d}<t_d^{-1/2}$.

Moreover, in the same way, we can prove that
$$
\|\chi^{d+1}\|_{k+2,r_d} \leq M_k
t_d^{-\frac{k-2s-4}{k-1}+A\frac{3s+5}{k-1}}\,,
$$
and we can assume (replacing $d_{k+1}$ by a higher integer if
necessary), that $\|\chi^{d+1}\|_{k+2,r_d}<t_d^{-1/2}$ for every
$d \geq d_{k+1}$.\\

$\bullet$ Proof of {\it (ii)} : Let $d\geq d_{k+1}$. Proceeding in
the same way as in the proof of Point  $(3_d)$ of the previous
lemma, we get
\begin{equation}
\|\Pi^{d+1}\|_{k+1} \leq
\|\Pi^d\|_{k+1}(1+M_k\|\chi^{d+1}\|_{k+2})^2
(1+p(\|\chi^{d+1}\|_{k+1}))\,,
\end{equation}
where $p$ is a polynomial with vanishing constant term. Now, since
$\|\chi^{d+1}\|_{k+1}$ and $\|\chi^{d+1}\|_{k+2}$ are strictly
smaller than $t_d^{-1/2}$, replacing $d_{k+1}$ by a higher integer
if necessary, we can assume that $\forall d \geq d_{k+1}$, we have
\begin{equation}
(1+M_k\|\chi^{d+1}\|_{k+2})^2
(1+p(\|\chi^{d+1}\|_{k+1}))<1+\frac{1}{(d+1)(d+2)}\,.
\end{equation}
We choose a positive constant ${\tilde C}_{k+1}$ such that
$\|\Pi^{d_{k+1}}\|_{k+1}<{\tilde C}_{k+1} \Big(
\frac{d_{k+1}+1}{d_{k+1}+2} \Big)$ and we can conclude by
induction, as in the previous lemma, that for all $d \geq
d_{k+1}$,
\begin{equation}
\|\Pi^d\|_{k+1,r_d}<{\tilde C}_{k+1} (1-\frac{1}{d+2})\,.
\end{equation}
Note that the constant ${\tilde C}_{k+1}$ is not the $C_{k+1}$ of
the lemma. Later, we will choose $C_{k+1}$ to be greater than
${\tilde C}_{k+1}$ and satisfying other conditions.\\

$\bullet$ Proof of  {\it (iii)} :  The idea is exactly the same as
in the previous lemma, using the interpolation inequality
(\ref{eqn:interpolation}) with $k$ and $2k-1$. Let $d\geq
d_{k+1}-1 \geq d_k$. By Lemmas \ref{lem:composition} and
\ref{lem:inverse}, in order to estimate
$\|f^{d+1}\|_{k+1,r_{d+1}}$ we just have to estimate
$\|\delta\big(\varphi^{d+1}+h(f^d)\big)+h(\delta
f^d)+Q^d\|_{k+1,r_d}$. As above, we write
\begin{eqnarray*}
\|\delta\big(\varphi^{d+1}+h(f^d)\big)\|_{k+1} &\leq& M_k\|h(f^d)-S(t_d)h(f^d)\|_{k+2}\\
 &\leq& M_k t_d^{-1}\|h(f^d)\|_{k+3}\quad {\mbox { by (\ref{eqn:smoothing2})}}\\
 &\leq& M_k t_d^{-1}\|f^d\|_{k+s+3}\quad {\mbox { by (\ref{eqn:estimate-h1})}}\\
 &\leq& M_k t_d^{-1}\|f^d\|_k^{\frac{k-s-4}{k-1}}\|f^d\|_{2k-1}^{\frac{s+3}{k-1}}\quad
{\mbox {by (\ref{eqn:interpolation})}}
\end{eqnarray*}
Then, since $\|f^d\|_{k}<C_k t_d^{-1}$ and $\|f^d\|_{2k-1}<C_k
t_d^{A}$ we have
$$
\|\delta\big(\varphi^{d+1}+h(f^d)\big)\|_{k+1} \leq M_k
t_d^{-1-\frac{k-s-4}{k-1}+ A\frac{s+3}{k-1}}\,.
$$
In the same way as in Point $(5_d)$ of the previous lemma, we can
estimate $\| h(\delta f^d) \|_{k+1}$ by $M_k
t_d^{-2\frac{k-s-3}{k-1}+2A\frac{s+2}{k-1}}$. Now, we just have to
estimate $\|Q^d\|_{k+1}$. By the definition of $Q^d$, we have
\begin{eqnarray*}
\|Q^d\|_{k+1} &\leq& M_k (\|f^d\|_{k+1}+\|g^d\|_{k+1}
+\|\Pi^d\|_{k+1}\|\varphi^{d+1}\|_{k+1}) \|\varphi^{d+1}\|_{k+2}\\
  &\leq& M_k (\|f^d\|_{k}^{\frac{k-2}{k-1}}\|f^d\|_{2k-1}^{\frac{1}{k-1}} +
\|g^d\|_{k}^{\frac{k-2}{k-1}}\|g^d\|_{2k-1}^{\frac{1}{k-1}} \\
  & & +  \|\Pi^d\|_{k+1}\|\varphi^{d+1}\|_{k+1}) \|\varphi^{d+1}\|_{k+2}
\end{eqnarray*}

We saw in {\it (ii)} that $\|\Pi^d\|_{k+1}\leq {\tilde C}_{k+1}$.
Moreover we saw in {\it (i)} that $\|\chi^{d+1}\|_{k+2,r_d} \leq
M_k t_d^{-\frac{k-2s-4}{k-1}+A\frac{3s+5}{k-1}}$. Since
$\frac{1}{k-1} < \frac{2s+3}{k-1}$, we obtain
$$
\|Q^d\|_{k+1} \leq M_k
t_d^{-2\frac{k-2s-4}{k-1}+2A\frac{3s+5}{k-1}}\,.
$$
We then obtain
$$
\|f^{d+1}\|_{k+1,r_{d+1}} \leq M_k
t_d^{-2\frac{k-2s-4}{k-1}+2A\frac{3s+5}{k-1}}\,,
$$
and we deduce, in the same way as in the proof of Point $(5_d)$ of
the previous lemma that for all $d\geq d_{k+1}-1$,
$$
\|f^{d+1}\|_{k+1,r_{d+1}} < M_k t_d^{-\mu}\,,
$$
where $-\mu<-3/2$. Therefore, replacing $d_{k+1}$ by a greater
integer if necessary, we have for every $d\geq d_{k+1}-1$
\begin{equation}
\|f^{d+1}\|_{k+1,r_{d+1}} < t_d^{-3/2}=t_{d+1}^{-1}\,.
\end{equation}
In the same way, we can show that
\begin{equation}
\|g^{d+1}\|_{k+1,r_{d+1}} < t_{d+1}^{-1}\,.
\end{equation}


$\bullet$ Proof of  {\it (iv)} : First recall that we have
\begin{eqnarray*}
f_{ij}^{d} &=& \{x_i,x_j\}_{d}-\sum_u c_{ij}^u x_u\\
g_{i\alpha}^{d} &=& \{x_i,y_\alpha\}_{d}-\sum_\beta
a_{i\alpha}^\beta y_\beta
\end{eqnarray*}
and, as in the proof of Point $(4_{d+1})$ of the previous lemma,
we can write
\begin{equation}
\begin{array}{l}
\|f^{d}\|_{2k+1,r_{d}} < V \|\Pi^{d}\|_{2k+1,r_{d}} \ , \cr
\|g^{d}\|_{2k+1,r_{d}} < V \|\Pi^{d}\|_{2k+1,r_{d}} \ ,
\end{array}
\label{eqn:est-fg2k+1}
\end{equation}
where $V>1$ is a positive constant independent of $d$ and $k$.\\

Now, we estimate $\|\Pi^{d+1}\|_{2k+1,r_{d+1}}$ for $d\geq
d_{k+1}$. Recall that we have
$$
{\{ x_i,x_j \}}_{d+1}={\{x_i+\varphi^{d+1}_i,x_j+ \varphi^{d+1}_j
\}}_d \circ \theta_{d+1}^{-1}
$$
and the same type of equality for  ${\{ x_i,y_\alpha \}}_{d+1}$
and ${\{y_\alpha,y_\beta\}}_{d+1}$.

Applying Lemmas \ref{lem:compositionL} and \ref{lem:inverse} we
obtain
\begin{eqnarray*}
\|{\{ x_i,x_j\}}_{d+1} \|_{2k+1,r_{d+1}} &\leq&
\|\{x_i+\varphi_i^{d+1},x_j+\varphi^{d+1}_j\}_d\|_{2k+1,r_d}
P_k(\|\chi^{d+1}\|_{k+1,r_d})\\
    &+&  \|\{x_i+\varphi_i^{d+1},x_j+\varphi^{d+1}_j\}_d\|_{k+1,r_d}
     \|\chi^{d+1}\|_{2k+1,r_d}Q_k(\|\chi^{d+1}\|_{k+1,r_d})\,,
\end{eqnarray*}
where $P_k$ and $Q_k$ are polynomials functions which do not
depend on $d$. In the same way as in the proof of $(2_d)$, since
$\|\Pi^d\|_{k+1,r_d}<{\tilde C}_{k+1} \frac{d+1}{d+2}$ and
$\|\chi^{d+1}\|_{k+2,r_d}< t_d^{-1/2}$, we can show that the term
$\|\{x_i+\varphi_i^{d+1},x_j+\varphi^{d+1}_j\}_d\|_{k+1}$ is
bounded. Therefore, we can write
$$
\|{\{ x_i,x_j\}}_{d+1} \|_{2k+1,r_{d+1}} \leq M_k\big(
\|\{x_i+\varphi_i^{d+1},x_j+\varphi^{d+1}_j\}_d\|_{2k+1,r_d} +
\|\chi^{d+1}\|_{2k+1,r_d}\big)\,.
$$
As in the proof of $(2_d)$ of the previous lemma, we first study
the term $\chi^{d+1}$. Actually, we will estimate
$\|\chi^{d+1}\|_{2k+2}$ rather than $\|\chi^{d+1}\|_{2k+1}$
because it will be used to estimate the terms of type
$\|\{x_i+\varphi_i^{d+1},x_j+\varphi^{d+1}_j\}_d\|_{2k+1}$. We
first write $\|\varphi^{d+1}\|_{2k+2}\leq M_k
t_d^{s+3}\|h(f^d)\|_{2k-1-s}$ by the property
(\ref{eqn:smoothing1}) of the smoothing operator. Using the
estimate of the homotopy operator (\ref{eqn:estimate-h1}), we
obtain $\|\varphi^{d+1}\|_{2k+2}\leq M_k t_d^{s+3}\|f^d\|_{2k-1}
\leq M_k t_d^{A+s+3}$. Now, we have
\begin{eqnarray*}
\|\psi^{d+1}\|_{2k+2} &\leq& M_k t_d^{3s+4}\|{\hat h}({\hat
g}^d)\|_{2k-3s-2}\,
{\mbox{  by }}\, (\ref{eqn:smoothing1}) \\
                    &\leq& M_k t_d^{3s+4}\|{\hat g}^d\|_{2k-2-2s} {\mbox{  by }}
\, (\ref{eqn:estimate-h1})
\end{eqnarray*}
Then the definition of ${\hat g}^d$, the Leibniz rule of
derivation and  (\ref{eqn:estimate-h1}) give
\begin{eqnarray*}
\|\psi^{d+1}\|_{2k+2} &\leq& M_k t_d^{3s+4}(\|g^d\|_{2k-2-2s}
+\|\Pi^d\|_{2k-2-2s}\|h(f^d)\|_{k-s-1+1}\\
                      & & +\|\Pi^d\|_{k-s-1}\|h(f^d)\|_{2k-2-2s+1})\\
                    &\leq& M_k t_d^{3s+4}(\|g^d\|_{2k-1}
+\|\Pi^d\|_{2k-1}\|f^d\|_{k}+\|\Pi^d\|_{k}\|f^d\|_{2k-1})\\
                    &\leq& M_k t_d^{A+3s+4}\,.
\end{eqnarray*}
Therefore, we can write
$$
\|\chi^{d+1}\|_{2k+2} \leq M_k t_d^{A+3s+4}\,.
$$

Now, in the same way as in the proof of Point $(2_d)$ of the
previous lemma, using the Leibniz formula of the derivation of a
product and the estimate $\|\chi^{d+1}\|_{k+2,r_d}<t_d^{-1/2}$, we
get
\begin{eqnarray*}
\|\{x_i+\varphi_i^{d+1},x_j+\varphi^{d+1}_j\}_d\|_{2k+1,r_d}
&\leq& M_k\big(\|\Pi^d\|_{2k+1}
(1+\|\varphi^{d+1}\|_{k+2})^2\\
& & \quad +\|\Pi^d\|_{k+1}(1+\|\varphi^{d+1}\|_{k+2})(1+\|\varphi^{d+1}\|_{2k+2})\big)\\
 &\leq& M_k(\|\Pi^d\|_{2k+1}+\|\varphi^{d+1}\|_{2k+2}+1)\\
 &\leq& M_k(\|\Pi^d\|_{2k+1}+t_d^{A+3s+4})\,.
\end{eqnarray*}
Consequently, we have
$$
\|{\{ x_i,x_j\}}_{d+1} \|_{2k+1,r_{d+1}} \leq M_k
(\|\Pi^d\|_{2k+1}+t_d^{A+3s+4})\,.
$$

In the same way, we can estimate
$\|\{x_i,y_\alpha\}_{d+1}\|_{2k+1,r_{d+1}}$ and
$\|\{y_\alpha,y_\beta\}_{d+1}\|_{2k+1,r_{d+1}}$ by
$M_k(\|\Pi^d\|_{2k+1}+t_d^{A+3s+4})$, which implies
$$
\|\Pi^{d+1}\|_{2k+1,r_{d+1}} \leq M_k
(\|\Pi^d\|_{2k+1,r_d}+t_d^{A+3s+4})\,.
$$

Finally, since $A>6s+8$, we can assume, replacing $d_{k+1}$ by a
higher integer if necessary, that $M_k t_d^{A+3s+4}
<\frac{1}{2V}t_d^{3A/2}$ for every $d \geq d_{k+1}$ (which also
implies that $M_k< \frac{1}{2V} t_d^{A/2}$). We then obtain,
$\forall d \geq d_{k+1}$,
\begin{equation}
\|\Pi^{d+1}\|_{2k+1,r_{d+1}} \leq \frac{1}{2V}
t_{d}^{A/2}\|\Pi^d\|_{2k+1,r_d} + \frac{1}{2V}t_{d}^{3A/2}\,.
\label{eqn:est-Pi2k+1}
\end{equation}

To conclude, if we choose a positive constant $C_{k+1}$ such that
$$
C_{k+1} > {\rm Max} \Big(1, {\tilde C}_{k+1},
\frac{\|\Pi^{d_{k+1}}\|_{2k+1,r_{d_{k+1}}}}{t_{d_{k+1}}^A}
\Big)\,,
$$
we then obtain, using (\ref{eqn:est-Pi2k+1}) and an induction,
$$
\|\Pi^{d}\|_{2k+1,r_{d}} < \frac{C_{k+1}}{V}t_{d}^{A} < C_{k+1}
t_{d}^{A}\,,
$$
for all $d\geq d_{k+1}$.

Finally, the estimates in (\ref{eqn:est-fg2k+1}) give, for all
$d\geq d_{k+1}$,
\begin{eqnarray*}
\|f^{d}\|_{2k+1,r_{d}} &<& C_{k+1}t_d^A\\
\|g^{d}\|_{2k+1,r_{d}} &<& C_{k+1}t_d^A\,.
\end{eqnarray*}

Moreover, the definition of $C_{k+1}$ completes the proof of the
points (i), (ii) and (iii).

Lemma \ref{lem:Hamilton2} is proved. \QED

\section{The case of Lie algebroids}
\label{section:Algebroids} In this section we briefly mention the
proof of Theorem \ref{thm:LeviAlgebroidS}. Similarly to the
analytic case (see \cite{Zung-Levi2002}), it is almost the same as
the proof of Theorem \ref{thm:LeviPoissonS}.

Let $A$ be a local $N$-dimensional smooth Lie algebroid (or
$C^{2q-1}$-smooth) over $(\bbR^n,0)$. We suppose that the anchor
map $\# : A \to T\bbR^n$, vanishes on $A_0$, the fiber of $A$ over
point $0$. It is well-known (see for instance
\cite{Cannas-Weinstein99}) that the Lie algebroid $A$ induces and
is, in fact, determined by a fiber-wise linear Poisson structure
on the dual bundle $A^\ast$.  More precisely, if
$(x_1,\hdots,x_n)$ is a local coordinate system on $\mathbb{R}^n$
and $(e_1,\hdots,e_N)$ is a local basis of sections, then
$(x_1,\hdots,x_n,e_1,\hdots,e_N)$ can be seen as a coordinate
system for $A^\ast$, which is linear on the fibers. The Poisson
structure on $A^\ast$ is given by
\begin{equation}
\begin{array}{l}
\{e_i,e_j\} = [e_i,e_j]\ , \cr \{e_i,x_j\} = \#e_i(x_j) \ ,\cr
\{x_i,x_j\} = 0\ .
\end{array}
\end{equation}
This Poisson structure is fiber-wise linear in the sense that the
bracket of two fiber-wise linear functions is again a fiber-wise
linear function, the bracket of a fiber-wise linear function with
a base function is a base function and the bracket of two base
functions is zero.

As in the statement of Theorem \ref{thm:LeviAlgebroidS}, we denote
by $\mathfrak{l}$ the $N$-dimensional Lie algebra in the linear
part of $A$ at $0$ (i.e. the isotropy algebra of $A$ at $0$), and
by $\mathfrak{g}$ a compact semisimple Lie subalgebra of
$\mathfrak{l}$. We can rewrite the basis of sections
$(e_1,\hdots,e_N)$ as $(s_1,\hdots,s_m,v_1,\hdots,s_{N-m})$ ($m$
is the dimension of $\mathfrak{g}$) where $(s_1,\hdots,s_m)$ span
$\mathfrak{g}$ and $(v_1,\hdots,v_{N-m})$ span a linear complement
of $\mathfrak{g}$ in $\mathfrak{l}$ which is invariant under the
adjoint action of $\mathfrak{g}$.

To prove Theorem \ref{thm:LeviAlgebroidS}, it suffices to find a
Levi factor for the dual Lie-Poisson structure which consists of
fiber-wise linear functions. The existence of a Levi factor is
given by Theorem \ref{thm:LeviPoissonS} and we only have to make
sure that this Levi factor can be chosen so that it consists of
fiber-wise linear functions. Actually the proof is the same as for
Theorem \ref{thm:LeviPoissonS} with few modifications :

The symbol $\mathcal{C}_r$ denotes now the subspace of the space
$C^\infty(B_r)$ of $C^\infty$-smooth real-valued functions on
$B_r$ (where $B_r \subset B_r^n \times \bbR^N$ is the closed ball
centered at $0$ and of radius $r$ in $\bbR^{n+N} = \bbR^n \times
\bbR^N$), which consists of fiber-wise linear functions vanishing
at $0$ whose first derivatives also vanish at $0$.

The symbol $\mathcal{Y}_r$ denotes now the space of
$C^\infty$-smooth vector fields on $B_r$ of the type
\begin{equation*}
\sum_{i=1}^{N-m} p_i \frac{\partial }{\partial v_i}+
\sum_{i=1}^{n} q_i \frac{\partial }{\partial x_i} \ ,
\end{equation*}
such that $p_i$ and $q_i$ vanish at $0$ and their first
derivatives also vanish at $0$ and, $p_i$ are fiber-wise linear
functions and $q_i$ are base functions.

One can check that these spaces are tame Fr\'echet spaces and
$\mathfrak{g}$-modules with the same actions as defined in Section
2. We then still have the homotopy operators and all the
properties we saw in Sections 2 and 3. The algorithm of
construction of the sequence of diffeomorphisms is the same as for
Theorem \ref{thm:LeviPoissonS} and one can check that if the
Poisson structure $\{ \,,\, \}_d$ is fiber-wise linear then $\{
\,,\, \}_{d+1}$ is fiber-wise linear too.

\section{Appendix: a Nash-Moser normal form theorem}

In this appendix we will generalize Theorem \ref{thm:LeviPoissonS}
into an abstract smooth normal form theorem, which we call a
\emph{Nash-Moser normal form theorem}, because its proof is
similar to the proof of Theorem \ref{thm:LeviPoissonS} and is
based on the Nash-Moser fast convergence method. Of course, Conn's
smooth linearization theorem \cite{Conn-Smooth1985}, as well as
our smooth Levi decomposition theorems, can be viewed as
particular cases of this abstract smooth normal form theorem,
modulo Lemma \ref{lem:HomotopyOp} about the norm of homotopy
operators. It is hoped that our abstract Nash-Moser normal form
theorem can be used or easily adapted for the study of other
smooth normal form problems as well.

\subsection{The setting}

Grosso modo, the situation is as follows: we have a group
${\mathcal G}$ (say of diffeomorphisms) which acts on a set $\cF$
(of structures). Inside $\cF$ there is a subset $\mathcal N$ (of
structures in normal form). We want to show that, under some
appropriate conditions, each structure can be put into normal
form, i.e. for each element $f \in \cF$ there is an element $\Phi
\in {\mathcal G}$ such that $\Phi.f \in {\mathcal N}$. We will
assume that $\cF$ is a subset of a linear space $\mathcal H$ (a
space of tensors) on which $\mathcal G$ acts, and $\mathcal N$ is
the intersection of $\cF$ with a linear subspace $\mathcal V$ of
$\mathcal H$. To formalize the situation involving smooth
\emph{local} structures (defined in a neighborhood of something),
let us introduce the following notions of \emph{SCI-spaces} and
\emph{SCI-groups}. Here SCI stands for \emph{scaled $C^\infty$
type}. Our aim here is not to create a very general setting, but
just a setting which works and which can hopefully be adjusted to
various situations. So our definitions below (especially the
inequalities appearing in them) are probably not ``optimal'', and
can be improved, relaxed, etc.

\noindent {\bf SCI-spaces.} An \emph{SCI-space} $\mathcal{H}$ is a
collection of Banach spaces
$(\mathcal{H}_{k,\rho},\|\,\|_{k,\rho})$ with $0<\rho\leq 1$ and
$k\in\Z_+ = \{0, 1, 2, \dots\}$ ($\rho$ is called the
\emph{radius} parameter, $k$ is called the \emph{smoothness
parameter}; we say that $f\in\mathcal{H}$ if
$f\in\mathcal{H}_{k,\rho}$ for some $k$ and $\rho$, and in that
case we say that $f$ is $k$-smooth and defined in radius $\rho$)
which satisfies the following properties:
\begin{itemize}
    \item If $k<k^\prime$, then for any $0 < \rho \leq 1$, $\mathcal{H}_{k^\prime,\rho}$ is a
linear subspace of $\mathcal{H}_{k,\rho}$:
$\mathcal{H}_{k^\prime,\rho} \subset \mathcal{H}_{k,\rho}$.
    \item If $0 < \rho<\rho^\prime \leq 1$, then for each $k \in \bbZ_+$, there
    is a given linear map, called the \emph{projection map}, or \emph{radius restriction
    map},
    $$\pi_{\rho,\rho'}: \mathcal{H}_{k,\rho^\prime} \rightarrow
    \mathcal{H}_{k,\rho}. $$
    These projections don't depend on $k$ and satisfy the natural
    commutativity condition $\pi_{\rho,\rho''} = \pi_{\rho',\rho''} \circ
    \pi_{\rho,\rho'}$. If $f \in \mathcal{H}_{k,\rho}$ and $\rho' < \rho$,
    then by abuse of language we will still denote by $f$ its projection to
    $\mathcal{H}_{k,\rho'}$ (when this notation does not lead to confusions).
    \item For any $f$ in $\mathcal H$ we have
    \begin{equation}
    \|f\|_{k,\rho}\geq\|f\|_{k^\prime,\rho^\prime} \ \ \forall \ k\geq k^\prime, \rho \geq \rho^\prime.
    \end{equation}
    In the above inequality, if $f$ is not in $\mathcal{H}_{k,\rho}$ then we
    put $\|f\|_{k,\rho} = + \infty$, and if $f$ is in
    $\mathcal{H}_{k,\rho}$ then the right hand side means the norm of the
    projection of $f$ to $\mathcal{H}_{k',\rho'}$, of course.
    \item There is a smoothing operator for each $\rho$, which depends
    continuously on $\rho$. More precisely, for each $0 < \rho \leq 1$ and
    each $t > 1$ there is a linear map, called the \emph{smoothing
    operator},
\begin{equation}
S_\rho(t):\mathcal{H}_{0,\rho} \longrightarrow
\mathcal{H}_{\infty,\rho}=\bigcap_{k=0}^\infty
\mathcal{H}_{k,\rho}
\end{equation}
which satisfies the following inequalities: for any $p, q \in
\bbZ_+$, $p \geq q$ we have
\begin{eqnarray}
\|S_\rho(t)f\|_{p,\rho} &\leq& C_{\rho,p,q}
t^{p-q}\|f\|_{q,\rho} \label{axiom:smoothing1}\\
\|f-S_\rho(t)f\|_{q,\rho} &\leq& C_{\rho,p,q}
t^{q-p}\|f\|_{p,\rho}\label{axiom:smoothing2}
\end{eqnarray}
where $C_{\rho,p,q}$ is a positive constant (which does not depend
on $f$ nor on $t$) and which is continuous with respect to $\rho$.
\end{itemize}

In the same way as for the Fr\'echet spaces (see for instance
\cite{Sergeraert1972}), the two properties
(\ref{axiom:smoothing1}) and (\ref{axiom:smoothing2}) of the
smoothing operator imply the following inequality called {\it
interpolation inequality} : for any positive integers $p$, $q$ and
$r$ with $p \geq q \geq r$ we have
\begin{equation}
(\|f\|_{q,\rho})^{p-r} \leq C_{p,q,r} (\|f\|_{r,\rho})^{p-q}
(\|f\|_{p,\rho})^{q-r}\,, \label{eqn:interpolNM}
\end{equation}
where $C_{p,q,r}$ is a positive constant which is continuous with
respect to $\rho$ and does not depend on $f$.\\

 Of course, if $\mathcal H$ is an SCI-space then each
$\mathcal{H}_{\infty,\rho}$ is a tame Fréchet space. The main
example that we have in mind is the space of functions in a
neighborhood of $0$ in the Euclidean space $\bbR^n$: here $\rho$
is the radius and $k$ is the smoothness class, i.e.
$\mathcal{H}_{k,\rho}$ is the space of $C^k$-functions on the
closed ball of radius $\rho$ and centered at $0$ in $\bbR^n$,
together with the maximal norm (of each function and its partial
derivatives up to order $k$); the projections are restrictions of
functions to balls of smaller radii.

By an \emph{SCI-subspace} of an SCI-space $\mathcal{H}$, we mean a
collection $\mathcal{V}$ of subspaces $\mathcal{V}_{k,\rho}$ of
$\mathcal{H}_{k,\rho}$, which themselves form an SCI-space (under
the induced norms, induced smoothing operators, induced inclusion
and radius restriction operators from $\mathcal{H}$ - it is
understood that these structural operators preserve
$\mathcal{V}$).

By a \emph{subset}  of an SCI-space $\mathcal{H}$, we mean a
collection $\mathcal{F}$ of subsets $\mathcal{F}_{k,\rho}$ of
$\mathcal{H}_{k,\rho}$, which are invariant under the inclusion
and radius restriction maps of $\mathcal{H}$.

\emph{Remark}. The above notion of SCI-spaces generalizes at the
same time the notion of tame Fréchet spaces and the notion of
scales of Banach spaces \cite{Zehnder1975}. Evidently, the scale
parameter is introduced to treat local problems. When things are
globally defined (say on a compact manifold), then the scale
parameter is not needed, i.e. $\mathcal{H}_{k,\rho}$ does not
depend on $\rho$ and we get back to the situation of tame Fréchet
spaces, as studied by Sergeraert \cite{Sergeraert1972} and
Hamilton \cite{Hamilton-Complex1977,Hamilton-NashMoser1982}.

\noindent {\bf SCI-groups.} An {\emph {SCI-group}} $\mathcal{G}$
consists of elements $\Phi$ which are written as a (formal) sum
\begin{equation}
\Phi = Id + \chi,
\end{equation}
where $\chi$ belongs to an SCI-space $\mathcal{W}$, together with
\emph{scaled group laws} to be made more precise below. We will
say that $\mathcal{G}$ is modelled on $\mathcal{W}$, if $\chi \in
\mathcal{W}_{k,\rho}$ then we say that $\Phi = Id + \chi \in
\mathcal{G}_{k,\rho}$ and $\chi = \Phi - Id$ (so as a space,
$\mathcal G$ is the same as $\mathcal W$, but shifted by $Id$),
$Id = Id + 0$ is the neutral element of $\mathcal{G}$.

\emph{Scaled composition (product) law}. There is a positive
constant $C$ (which does not depend on $\rho$ or $k$) such that if
$0 < \rho' < \rho \leq 1$, $k \geq 1$, and $\Phi=Id+\chi \in
\mathcal{G}_{k,\rho}$ and $\Psi=Id+\xi \in \mathcal{G}_{k,\rho}$
such that
\begin{equation}
\rho^\prime/\rho\leq 1-C\|\xi\|_{1,\rho} \label{cond:radius}
\end{equation}
then we can compose $\Phi$ and $\Psi$ to get an element $\Phi
\circ \Psi$ with $\|\Phi\circ\Psi-Id\|_{k,\rho^\prime}<\infty$,
i.e. $\Phi\circ\Psi$ can be considered as an element of
$\mathcal{G}_{k,\rho'}$ (if $\rho'' < \rho'$ then of course
$\Phi\circ\Psi$ can also be considered as an element of
$\mathcal{G}_{k,\rho''}$, by the restriction of radius from
$\rho'$ to $\rho''$). Of course, we require the composition to be
\emph{associative} (after appropriate restrictions of radii).


\emph{Scaled inversion law}. There is a positive constant $C$ (for
simplicity, take it to be the same constant as in Inequality
(\ref{cond:radius})) such that if $\Phi \in \mathcal{G}_{k,\rho}$
such that
\begin{equation} \label{cond:radius2}
\|\Phi - Id\|_{1,\rho} < 1/C
\end{equation}
then we can define an element, denoted by $\Phi^{-1}$ and called
the inversion of $\Phi$, in $\mathcal{G}_{k,\rho'}$, where $\rho'
= (1- {1 \over 2}C\|\Phi - Id\|_{1,\rho}) \rho $, which satisfies
the following condition: the compositions $\Phi \circ \Phi^{-1}$
and $\Phi^{-1} \circ \Phi$ are well-defined in radius $\rho'' =
(1- C\|\Phi - Id\|_{1,\rho})\rho$ and coincide with the neutral
element $Id$ there.

\emph{Continuity conditions}. We require that the above scaled
group laws satisfy the following continuity conditions i), ii) and
iii) in order for $\mathcal{G}$ to be called an SCI-group.

i) For each $k \geq 1$ there is a polynomial $P = P_k$ (of one
variable), such that for any $\chi \in \mathcal{W}_{2k-1,\rho}$
with $\|\chi\|_{1,\rho} < 1/C$ we have
\begin{equation}
\|(Id+\chi)^{-1}-Id\|_{k,\rho^\prime}\leq \|\chi\|_{k,\rho}
P(\|\chi\|_{k,\rho}) \label{axiom:inverse}\ ,
\end{equation}
where $\rho' = (1 - C\|\chi\|_{1,\rho})\rho$.

ii) If $(\Phi_m)_{m\geq 0}$ is a sequence in
$\mathcal{G}_{k,\rho}$ which converges (with respect to
$\|\,\|_{k,\rho}$) to $\Phi$, then the sequence
$(\Phi_m^{-1})_{m\geq 0}$ also converges to $\Phi^{-1}$ in
$\mathcal{G}_{k,\rho^\prime}$, where $\rho' = (1 - C\|\Phi -
Id\|_{1,\rho})\rho$.


iii) For each $k \geq 1$ there are polynomials $P$ and $Q$ (of one
variable) with vanishing constant term such that if $\Phi=Id+\chi$
and $\Psi = Id + \xi$ are in $\mathcal{G}_{k,\rho}$ and if
$\rho^\prime$ and $\rho$ satisfy Relation (\ref{cond:radius}),
then we have
\begin{equation}
\|\Phi\circ\Psi-\Phi\|_{k,\rho^\prime}\leq P(\|\xi\|_{k,\rho})+
\|\chi\|_{k+1,\rho} Q(\|\xi\|_{k,\rho})\,. \label{axiom:product}
\end{equation}

{\it Remark :} As a consequence of the last condition we have,
with the same notations, the following inequality:
\begin{equation}
\|\Phi\circ\Psi-Id\|_{k,\rho^\prime}\leq P(\|\xi\|_{k,\rho})+
\|\chi\|_{k+1,\rho}(1+Q(\|\xi\|_{k,\rho}))\,.
\label{conseq-axiom-product}
\end{equation}

\noindent{\bf SCI-actions}. We will say that there is a
\emph{linear left SCI-action} of an SCI-group $\mathcal{G}$ on an
SCI-space $\mathcal{H}$ if there is a positive integer $\gamma$
(and a positive constant $C$) such that, for each $\Phi = Id +
\chi \in \mathcal{G}_{k,\rho}$ and $f \in \mathcal{H}_{k,\rho'}$
with $\rho' = (1 - C\|\chi\|_{1,\rho})\rho$, the element $\Phi. f$
(the image of the action of $\Phi$ on $f$) is well-defined in
$\mathcal{H}_{k,\rho'}$, the usual axioms of a left group action
modulo appropriate restrictions of radii (so we have \emph{scaled
action laws}) are satisfied, and the following three inequalities
i), ii), iii) expressing some continuity conditions are also
satisfied:

i) For each $k$ there is a polynomial function $P = P_k$ with
vanishing constant term such that
\begin{equation} \label{axiom:action3}
\|(Id+\chi)\cdot f\|_{k,\rho^\prime} \leq \|f\|_{k,\rho} (1+
P(\|\chi\|_{k+\gamma,\rho})) \ .
\end{equation}

ii) For each $k$ there are polynomials $Q$ and $R$ (which depend
on $k$) such that
\begin{equation}
\|(Id+\chi)\cdot f\|_{2k-1,\rho^\prime} \leq \|f\|_{2k-1,\rho}
Q(\|\chi\|_{k+\gamma,\rho})+\|\chi\|_{2k-1+\gamma,\rho}
\|f\|_{k,\rho} R(\|\chi\|_{k+\gamma,\rho}) \label{axiom:action2}
\end{equation}

iii) There is a polynomial function $O$ of 2 variables such that
\begin{equation}
\|(\Phi+\chi)\cdot f-\Phi\cdot f\|_{k,\rho^\prime} \leq
\|\chi\|_{k+\gamma,\rho} \|f\|_{k+\gamma,\rho} O(\|\Phi -
Id\|_{k+\gamma,\rho}, \|\chi\|_{k+\gamma,\rho})
\label{axiom:action1}
\end{equation}

In the above inequalities, $\rho'$ is related to $\rho$ by a
formula of the type $\rho' = \left( 1 - C(\|\chi\|_{1,\rho}+
\|\Phi -Id\|_{1,\rho})\right) \rho$. ($\Phi = Id$ in the first two
inequalities).

The main example of a (linear left) SCI-action that we have in
mind is the push-forward action of the SCI-group of local
diffeomorphisms of $(\bbR^n, 0)$ on the SCI-space of local tensors
of a given type (e.g. 2-vector fields) on $(\bbR^n, 0)$.

\subsection{Normal form theorem} Roughly speaking, the following theorem
says that whenever we have a ``fast normalizing algorithm'' in an
SCI setting then it will lead to the existence of a smooth
normalization. ``Fast'' means that, putting loss of
differentiability aside, one can ``quadratize'' the error term at
each step (going from ``$\epsilon$-small'' error to
``$\epsilon^2$-small'' error).

In the statement of the following theorem, the polynomials 
$P_k$, $Q_k$, $R_k$ and $T_k$ depend on $k$ and may depend 
on $\rho$ continuously, but do not depend on $f$.

\begin{thm}
Let $\mathcal{H}$ be a SCI-space, $\mathcal{V}$ a SCI-subspace of
$\mathcal{H}$, and $\mathcal{F}$ a subset of $\mathcal{H}$,
$\mathcal{F} \ni 0$. Denote
$\mathcal{N}=\mathcal{F}\cap\mathcal{V}$. Assume that there is a
projection $\pi: \mathcal{H}\longrightarrow \mathcal{V}$
(compatible with restriction and inclusion maps) such that for
every $f$ in $\mathcal{H}_{k,\rho}$, the element
$\zeta(f)=f-\pi(f)$ satisfies 
\begin{equation}
\|\zeta(f)\|_{k,\rho}\leq
\|f\|_{k,\rho} T_k(\|f\|_{[(k+1)/2],\rho})
\label{eqn:proj}
\end{equation}
for all $k \in \bbN$ (or at least for all $k$ sufficiently large), where
$[\;]$ is the integer part and $T_k$ a polynomial. Let $\mathcal{G}$ be an
SCI-group acting on $\mathcal{H}$ by a linear left SCI-action
which preserves $\mathcal{F}$.  Assume that there is $s\in\N$  such
that for every $f$ in $\mathcal{F}$ and $0 < \rho \leq 1$, there
is an element $\Phi_f=Id+\chi_f\in\mathcal{G}$ (which may depend
on $\rho$ but doesn't depend on $k$) such that for all $k$ in
$\bbN$ (or at least for all $k$ sufficiently large), 
\begin{eqnarray}
\|\chi_f\|_{k,\rho} &\leq& \|\zeta(f)\|_{k+s,\rho}
P_k(\|f\|_{[(k+1)/2]+s,\rho}) \label{eqn:estimate-chif}\\
  &+&
\|f\|_{k+s,\rho}\|\zeta(f)\|_{[(k+1)/2]+s,\rho}Q_k(\|f\|_{[(k+1)/2]+s,\rho})
\ , \nonumber
\end{eqnarray}
and that the element $f' :=\Phi_f\cdot f \in \mathcal{F}$
satisfies the inequality
\begin{equation}
\|\zeta(f')\|_{k,\rho'} \leq \|\zeta(f)\|_{k+s,\rho}^2 
R_k(\|f\|_{k+s,\rho},\|\chi_f\|_{k+s,\rho},\|f\|_{k,\rho})
\label{eqn:estimate-zeta} 
\end{equation}
($\rho$ and $\rho'$ verify Relation (\ref{cond:radius})) where
$P_k$ and $Q_k$ (resp. $R_k$) are some
polynomials of 1 variable (resp. 3 variables) and
the degree in the first variable of the polynomial $R_k$ does not depend on $k$. 
Then there exist $l\in\N$ and two positive constants
$\alpha$ and $\beta$ with the following property: for all $p \in
\bbN \cup \{\infty\}, p \geq l$, and for all
$f\in\mathcal{F}_{2p-1,\rho}$ with $\|f\|_{2l-1,\rho}<\alpha$ and
$\|f-0\|_{l,\rho}<\beta$, there exists
$\Psi\in\mathcal{G}_{p,\rho/2}$ such that $\Psi\cdot f\in
\mathcal{N}_{p,\rho/2}$. \label{thm:Nash-Moser}
\end{thm}

\emph{Proof.} We construct, by induction, a sequence
${(\Psi_d)}_{d\geq 1}$ in $\mathcal{G}$, and then a sequence
$f^d:= \Psi_{d}\cdot f$ in $\mathcal{F}$, which converges to
$\Psi\in \mathcal{G}_{p,\rho/2}$ and such that
$f^\infty:=\Psi\cdot f\in \mathcal{N}_{p,\rho/2}$.

In order to simplify, we can assume that the constant $s$ of the
theorem is the same as the integer $\gamma$ defined by the
SCI-action of $\mathcal{G}$ on $\mathcal{H}$ (see
(\ref{axiom:action3}), (\ref{axiom:action3}) and
(\ref{axiom:action3})). We first fix some parameters. Let $A=6s+5$
(actually, $A$ just have to be strictly larger than $6s+4$). 
We denote by $\tau$ the degree in the first variable of the
polynomials $R_k$ introduced in Theorem \ref{thm:Nash-Moser}. We
consider a positive real number $\varepsilon <1$ such that
\begin{equation}
-(1-\varepsilon)+A(1+\frac{\tau}{2})\varepsilon<-\frac{3}{4}\,.
\label{cond:Aepsilon}
\end{equation}
Finally, we fix a positive integer $l>3s+3$ which satisfies
\begin{equation}
\frac{2s+2}{l-1}<\varepsilon\,. \label{cond:lsespsilon}
\end{equation}

The construction of the sequences is the following : Let $t_0>1$
be a real constant~; this constant is still not really fixed and
will be chosen according to Lemma \ref{lem:NM1}. We then define
the sequence ${(t_d)}_{d\geq 1}$ by $t_{d+1}:=t_d^{3/2}$. We also
define the sequence $r_d:=(1+\frac{1}{d+1})\rho/2$. This is a
decreasing sequence such that $\rho/2 \leq r_d\leq  \rho$ for all
$d$. Note that we have $r_{d+1}=r_d(1-\frac{1}{(d+2)^2})$.

Let $p\geq l$ and $f$ in $\mathcal{F}_{2p-1,\rho}$. We start with
$f_0:=f\in \mathcal{F}_{2p-1,\rho}$. Now, assume that we have
constructed $f^d\in \mathcal{F}_{2p-1,r_d}$ for $d\geq 0$. We put
$\Phi_{d+1}:=\Phi_{f^d}=Id+\chi^{d+1}$ and ${\hat
\Phi}_{d+1}:=S(t_d)\Phi_{d+1}=Id+{\hat \chi}^{d+1}$. Then,
$f^{d+1}$ is defined by
$$
f^{d+1}={\hat \Phi}_{d+1}\cdot f^d\,.
$$
Roughly speaking, the idea is that the sequence ${(f^d)}_{d\geq
0}$ will be such that
$$\|\zeta(f^{d+1})\|_{p,r_{d+1}} \leq
\|\zeta(f^d)\|_{p,r_d}^2\,.$$

For every $d\geq 1$, we put $\Psi_d={\hat
\Phi}_d\circ\hdots\circ{\hat \Phi}_1$. We then have to show that
we can choose two positive constants $\alpha$ and $\beta$ such
that if $\|f\|_{2l-1,\rho}\leq\alpha$ and
$\|f-0\|_{l,\rho}\leq\beta$ then, the sequence ${(\Psi_d)}_{d\geq
1}$ converges with respect to $\|\,\|_{p,\rho/2}$. It will follow
from these two technical lemmas that we will prove later :

\begin{lem}
There exists a real number $t_0>1$ such that for any
$f\in\mathcal{F}_{2p-1,\rho}$ satisfying the conditions
$\|f^0\|_{2l-1,r_0}< t_0^A$, $\|\zeta(f^0)\|_{2l-1,r_0}< t_0^A$
and $\|\zeta(f^0)\|_{l,r_0}< t_0^{-1}$ then, with the construction
above, we have for all $d\geq 0$,
\begin{itemize}
    \item [$(1_d)$] $\quad \|{\hat \chi}^{d+1}\|_{l+s,r_d}< t_d^{-1/2}$
    \item [$(2_d)$] $\quad \|f^d\|_{l,r_d}< C \frac{d+1}{d+2}$ where $C$ is a
    positive constant
    \item [$(3_d)$] $\quad \|f^d\|_{2l-1,r_d} < t_d^A$
    \item [$(4_d)$] $\quad \|\zeta(f^d)\|_{2l-1,r_d}< t_d^A$
    \item [$(5_d)$] $\quad \|\zeta(f^d)\|_{l,r_d}< t_d^{-1}$
\end{itemize}
\label{lem:NM1}
\end{lem}

\begin{lem}
Suppose that for an integer $k\geq l$, there exists a constant
$C_k$ and an integer $d_k\geq 0$ such that for any $d\geq d_k$ we
have $\|f^d\|_{2k-1,r_d}<C_k t_d^A$,
$\|\zeta(f^d)\|_{2k-1,r_d}<C_k t_d^A$, $\|f^d\|_{k,r_d}<
C_k\frac{d+1}{d+2}$ and $\|\zeta(f^d)\|_{k,r_d}< C_k t_d^{-1}$.
Then, there exists a positive constant $C_{k+1}$ and an integer
$d_{k+1} > d_k$ such that for any $d\geq d_{k+1}$ we have
\begin{itemize}
    \item [(i)]   $\quad \|{\hat \chi}^{d+1}\|_{k+1+s,r_d}< C_{k+1} t_d^{-1/2}$
    \item [(ii)]  $\quad \|f^d\|_{k+1,r_d}< C_{k+1}\frac{d+1}{d+2}$
    \item [(iii)] $\quad \|f^d\|_{2k+1,r_d} < C_{k+1} t_d^A$
    \item [(iv)]  $\quad \|\zeta(f^d)\|_{2k+1,r_d}< C_{k+1} t_d^A$
    \item [(v)]   $\quad \|\zeta(f^d)\|_{k+1,r_d}< C_{k+1} t_d^{-1}$
\end{itemize}
\label{lem:NM2}
\end{lem}

{\it End of the proof of Theorem \ref{thm:Nash-Moser} : } We
choose $t_0$ as in Lemma \ref{lem:NM1}. Then we fix two positive
constants $\alpha$ and $\beta$ such that $t_0^A\geq \alpha$ and
$t_0^{-1}\geq \beta$. Now, if $f\in\mathcal{F}_{2p-1,\rho}$
satisfies $\|f\|_{2l-1,\rho}\leq\alpha$ and
$\|f-0\|_{l,\rho}\leq\beta$ then, since $\|\zeta(f)\|_{l,\rho}\leq
\|f-0\|_{l,\rho}$, using Lemma \ref{lem:NM1} and then applying
Lemma \ref{lem:NM2} repetitively, there exists a positive integer
$d_p$ such that for all $d\geq d_p$,
$$
\|{\hat \chi}^{d+1}\|_{p,r_d}< C_p t_d^{-1/2}\,.
$$

Actually it is more convenient to prove the convergence of the
sequence ${(\Psi_d^{-1})}_{d\geq 1}$. The point ii) of the
continuity conditions in the definition of SCI-group will then
give the convergence of ${(\Psi_d)}_{d\geq 1}$. For all positive
integer $d$, we have $\Psi_d^{-1}={\hat
\Phi}_1^{-1}\circ\hdots\circ{\hat \Phi}_d^{-1}$ and if we denote
${\hat \Phi}_d^{-1}=Id+{\hat \xi}^{d}$, the axiom
(\ref{axiom:inverse}) implies
$$
\|{\hat \xi}^{d+1}\|_{p,r_d}< M_p t_d^{-1/2}\,,
$$
for all $d\geq d_p$, where $M_p$ is a positive constant
independent of $d$. Now, by the inequality
(\ref{conseq-axiom-product}), the sequence
${(\Psi_d^{-1}-Id)}_{d\geq 1}$ is bounded and
(\ref{axiom:product}) gives then the
$\|\,\|_{p,\rho/2}$-convergence of ${(\Psi_d^{-1})}_{d\geq 1}$.


{\it Proof of Lemma \ref{lem:NM1} : } In this proof $M$ denotes a
positive constant which does not depend on $d$ and which varies
from inequality to inequality. As in the case of Poisson
structures, we prove this lemma by induction.

At the step $d=0$   
the only thing we have to verify is the point $(1_0)$ (for the
point $(3_0)$ we just choose the constant $C$ such that
$C>2\|f^0\|_{l,r_0}$) . We have $\|{\hat
\chi}^1\|_{l+s,r_0}=\|S(t_0) \chi^1\|_{l+s,r_0}\leq M
\|\chi^1\|_{l+s,r_0}$ by (\ref{axiom:smoothing1}). Therefore,
using (\ref{eqn:estimate-chif}) with the relation $l>3s+3$, and
the interpolation inequality (\ref{eqn:interpolNM}), we obtain
($P$ and $Q$ are polynomial functions) :
\begin{eqnarray*}
\|{\hat \chi}^1\|_{l+s,r_0} &\leq& M \|\zeta(f^0)\|_{l+2s,r_0}
P(\|f^0\|_{l,r_0})\\
   &+& M \|f^0\|_{l+2s,r_0}\|\zeta(f^0)\|_{l,r_0}
Q(\|f^0\|_{l,r_0})\\
   &\leq& M
\|\zeta(f^0)\|_{l,r_0}^{\frac{l-2s-1}{l-1}}
\|\zeta(f^0)\|_{2l-1,r_0}^{\frac{2s}{l-1}} \\
   &+& M
\|f^0\|_{l,r_0}^{\frac{l-2s-1}{l-1}}\|f^0\|_{2l-1,r_0}^{\frac{2s}{l-1}}
\|\zeta(f^0)\|_{l,r_0}\\
                          &\leq& M (t_0^{-\frac{l-2s-1}{l-1}+A\frac{2s}{l-1}}
                          +
                          t_0^{-1+A\frac{2s}{l-1}})
\end{eqnarray*}
Then, by (\ref{cond:lsespsilon}) and (\ref{cond:Aepsilon}), we
obtain $\|{\hat \chi}^1\|_{l+s,r_0} \leq M t_0^{-\mu}$ with
$-\mu<-3/4<-1/2$ and, replacing $t_0$ by a larger number if
necessary (independently of $f$ and $d$), we have $\|{\hat
\chi}^1\|_{l+s,r_0}<t_0^{-1/2}$. Note that we also proved that
$\|\chi^1\|_{l+s,r_0}<t_0^{-1/2}$.\\

Now, we suppose that the conditions $(1_d)\hdots (5_d)$ are
satisfied and we study the step $d+1$. The point $(1_{d+1})$ can
be proved as above.\\

Proof of $(2_{d+1})$ : According to (\ref{axiom:action3}) we have
$\|f^{d+1}\|_{l,r_{d+1}}\leq \|f^d\|_{l,r_d}(1+P(\|{\hat
\chi}^{d+1}\|_{l+s,r_d}))$ where $P$ is a polynomial with
vanishing constant term. Since $\|{\hat
\chi}^{d+1}\|_{l+s,r_d}<t_d^{-1}$ we can assume, choosing $t_0$
large enough, that $P(\|{\hat \chi}^{d+1}\|_{l+s,r_d})\leq
\frac{1}{(d+1)(d+3)}$ and we get
$$
\|f^{d+1}\|_{l,r_{d+1}}<C\frac{d+1}{d+2}(1+\frac{1}{(d+1)(d+3)})<
C\frac{d+2}{d+3}\,.
$$

Proof of $(3_{d+1})$ : We have $f^{d+1}={\hat {\Phi}}_{d+1}\cdot
f^d$ with ${\hat {\Phi}}_{d+1}=Id+{\hat
\chi}^{d+1}=Id+S(t_d)\chi^{d+1}$ thus, (\ref{axiom:action2}) gives
$$
\|f^{d+1}\|_{2l-1,r_{d+1}} \leq \|f^d\|_{2l-1,r_d} P_1(\|{\hat
{\chi}}^{d+1}\|_{l+s,r_d})+\|{\hat {\chi}}^{d+1}\|_{2l-1+s,r_d}
\|f^d\|_{l,r_d} P_2(\|{\hat \chi}^{d+1}\|_{l+s,r_d})
$$
where $P_1$ and $P_2$ are two polynomials. This gives, by $(1_d)$
and $(2_d)$,
$$
\|f^{d+1}\|_{2l-1,r_{d+1}} \leq M ( \|f^d\|_{2l-1,r_d} + \|{\hat
{\chi}}^{d+1}\|_{2l-1+s,r_d} )\,.
$$
Now, we have
\begin{eqnarray*}
\|{\hat {\chi}}^{d+1}\|_{2l-1+s,r_d} &\leq& Mt_d^{3s}
\|\chi^{d+1}\|_{2l-1-2s,r_d}\quad {\mbox { by }}\,
(\ref{axiom:smoothing1})\\
   &\leq& Mt_d^{3s} (\|\zeta(f^d)\|_{2l-1-s,r_d}
   P_3(\|f^d\|_{l,r_d})\\
   & & +\|f^d\|_{2l-1-s,r_d}\|\zeta(f^d)\|_{l,r_d}
   P_4(\|f^d\|_{l,r_d}))
   \quad {\mbox { by }}\,(\ref{eqn:estimate-chif})
\end{eqnarray*}
where $P_3$ and $P_4$ are polynomials. We get $\|{\hat
{\chi}}^{d+1}\|_{2l-1+s,r_d} \leq M t_d^{A+3s}$ and, consequently,
$$
\|f^{d+1}\|_{2l-1,r_{d+1}} \leq M t_d^{A+3s}\,.
$$
To finish, since $A=6s+5$, we have that
$\|f^{d+1}\|_{2l-1,r_{d+1}} \leq M t_d^{B}$ with $0<B<3A/2$ thus,
replacing $t_0$ by a larger number if necessary, we get
$\|f^{d+1}\|_{2l-1,r_{d+1}} < t_d^{3A/2}=t_{d+1}^A$.\\

Proof of $(4_{d+1})$ : We have
$$
\|\zeta(f^{d+1})\|_{2l-1,r_{d+1}}\leq
\|f^{d+1}\|_{2l-1,r_{d+1}} T(\|f^{d+1}\|_{l,r_{d+1}})
$$
where $T$ is a polynomial (see (\ref{eqn:proj})). Using the estimate of $(3_{d+1})$ and $(2_{d+1})$,
we obtain $\|\zeta(f^{d+1})\|_{2l-1,r_{d+1}}\leq M t_d^{A+3s}$,
and we conclude as above.

Proof of $(5_{d+1})$ : Recall that we have
$\Phi_{d+1}=Id+\chi^{d+1}$ and ${\hat
\Phi}_{d+1}=Id+S(t_d)\chi^{d+1}$. We can write
\begin{eqnarray*}
\|\zeta(f^{d+1})\|_{l,r_{d+1}} &=& \|\zeta({\hat \Phi}_{d+1}\cdot
f^{d})\|_{l,r_{d+1}}\\
 &\leq& \|\zeta({\hat \Phi}_{d+1}\cdot
f^{d}-\Phi_{d+1}\cdot f^{d})\|_{l,r_{d+1}}
+\|\zeta(\Phi_{d+1}\cdot f^{d})\|_{l,r_{d+1}}
\end{eqnarray*}
On the one hand, by (\ref{eqn:estimate-zeta}) and using the interpolation inequality
(\ref{eqn:interpolNM}), Point $(2_d)$
and the estimate $\|\chi^{d+1}\|_{l+s,r_d}<t_d^{-1/2}$ (see the
proof of $(1_0)$), we have
\begin{eqnarray*}
\|\zeta(\Phi_{d+1}\cdot f^{d})\|_{l,r_{d+1}} &\leq&
\|\zeta(f^d)\|_{l+s,r_d}^2 R_l(\|f^d\|_{l+s,r_d},\|\chi^{d+1}\|_{l+s,r_d},\|f^d\|_{l,r_d})\\
   &\leq& M\|\zeta(f^d)\|_{l,r_d}^{2\frac{l-s-1}{l-1}}
          \|\zeta(f^d)\|_{2l-1,r_d}^{2\frac{s}{l-1}} \times \\
  & &  R_l(\|f^d\|_{l,r_d}^{\frac{l-s-1}{l-1}}
\|f^d\|_{2l-1,r_d}^{\frac{s}{l-1}},\|\chi^{d+1}\|_{l+s,r_d},\|f^d\|_{l,r_d}) \\
 &\leq& M t_d^{-2\frac{l-s-1}{l-1}+2A\frac{s}{l-1}+A\frac{\tau s}{l-1}}
\end{eqnarray*}
recall that $\tau$ is the degree in the first variable of $R_l$.
Then, by (\ref{cond:lsespsilon}) and (\ref{cond:Aepsilon}), we
have $\|\zeta(\Phi_{d+1}\cdot f^{d})\|_{l,r_{d+1}}\leq
Mt_d^{-\mu}$ where $-\mu<-3/2$ and, replacing $t_0$ by a larger
number if necessary, we have $\|\zeta(\Phi_{d+1}\cdot
f^{d})\|_{l,r_{d+1}} <\frac{1}{2}t_d^{-3/2}$.

On the other hand, by (\ref{eqn:proj}),
\begin{eqnarray*}
\|\zeta({\hat \Phi}_{d+1}\cdot f^{d}-\Phi_{d+1}\cdot
f^{d})\|_{l,r_{d+1}} & \leq & \|{\hat \Phi}_{d+1}\cdot
f^{d}-\Phi_{d+1}\cdot f^{d}\|_{l,r_{d+1}} \times \\
  &  &  \times T(\|{\hat \Phi}_{d+1}\cdot
f^{d}-\Phi_{d+1}\cdot f^{d}\|_{l,r_{d+1}})
\end{eqnarray*}
and since ${\hat \Phi}_{d+1}=\Phi_{d+1}+({\hat
\chi}^{d+1}-\chi^{d+1})$, we have by (\ref{axiom:action1}),
\begin{eqnarray*}
\|{\hat \Phi}_{d+1}\cdot
f^{d}-\Phi_{d+1}\cdot f^{d}\|_{l,r_{d+1}} &\leq& \|{\hat
 \chi}^{d+1}-\chi^{d+1}\|_{l+s,r_d}\|f^d\|_{l+s,r_d}\\
 & & O(\|\chi^{d+1}\|_{l+s,r_d},\|{\hat \chi}^{d+1}-\chi^{d+1}\|_{l+s,r_d})
\end{eqnarray*}
where $O$ is a polynomial of 2 variables. Since $\|{\hat
\chi}^{d+1}\|_{l+s,r_d}$ and $\|\chi^{d+1}\|_{l+s,r_d}$ are both
majored by $t_d^{-1/2}$, we can write
$$
\|{\hat \Phi}_{d+1}\cdot f^{d}-\Phi_{d+1}\cdot
f^{d}\|_{l,r_{d+1}}\leq M \|{\hat
 \chi}^{d+1}-\chi^{d+1}\|_{l+s,r_d}\|f^d\|_{l+s,r_d}\,.
$$
By the interpolation inequality we can write
$\|f^d\|_{l+s,r_d}\leq M \|f^d\|_{2l-1,r_d}^{\frac{s}{l-1}}$.
Moreover, using the property (\ref{axiom:smoothing2}), the
estimate (\ref{eqn:estimate-chif}) with the inequality $l>3s+3$,
and then the interpolation inequality (\ref{eqn:interpolNM}), we
get
\begin{eqnarray*}
\|{\hat \chi}^{d+1}-\chi^{d+1}\|_{l+s,r_d} &\leq& M t_d^{-1}
\|\chi^{d+1}\|_{l+s+1,r_d}\\
 &\leq& M t_d^{-1} (\|\zeta(f^d)\|_{l+2s+1,r_d}
 P(\|f^d\|_{l,r_d})\\
 & & + \|f^d\|_{l+2s+1,r_d}\|\zeta(f^d)\|_{l,r_d} Q(\|f^d\|_{l,r_d}))\\
 &\leq& M t_d^{-1} (\|\zeta(f^d)\|_{l,r_d}^{\frac{l-2s-2}{l-1}}
\|\zeta(f^d)\|_{2l-1,r_d}^{\frac{2s+1}{l-1}}\\
 & & + \|f^d\|_{l,r_d}^{\frac{l-2s-2}{l-1}}
\|f^d\|_{2l-1,r_d}^{\frac{2s+1}{l-1}}\|\zeta(f^d)\|_{l,r_d})\\
 &\leq& M t_d^{-1} (t_d^{-\frac{l-2s-2}{l-1}+A\frac{2s+1}{l-1}}
+t_d^{-1+A\frac{2s+1}{l-1}})\,.
\end{eqnarray*}

Consequently, we have
\begin{eqnarray*}
\|{\hat \Phi}_{d+1}\cdot f^{d}-\Phi_{d+1}\cdot
f^{d}\|_{l,r_{d+1}} &\leq& M t_d^{A\frac{s}{l-1}} \|{\hat
\chi}^{d+1}-\chi^{d+1}\|_{l+s,r_d}\\
 &\leq& M t_d^{-2+A\frac{3s+1}{l-1}}\,,
\end{eqnarray*}
which implies
$$
\|\zeta({\hat \Phi}_{d+1}\cdot f^{d}-\Phi_{d+1}\cdot
f^{d})\|_{l,r_{d+1}}\leq  M t_d^{-2+A\frac{3s+1}{l-1}} 
T(M t_d^{-2+A\frac{3s+1}{l-1}})\,.
$$
As above, we can conclude that
$$
\|\zeta({\hat \Phi}_{d+1}\cdot f^{d+1}-\Phi_{d+1}\cdot
f^{d+1})\|_{l,r_{d+1}}< \frac{1}{2} t_d^{-3/2}\,.
$$
Finally, we obtain
\begin{eqnarray*}
\|\zeta(f^{d+1})\|_{l,r_{d+1}} &\leq& \|\zeta({\hat
\Phi}_{d+1}\cdot f^{d}-\Phi_{d+1}\cdot f^{d})\|_{l,r_{d+1}}
+\|\zeta(\Phi_{d+1}\cdot f^{d})\|_{l,r_{d+1}}\\
   &<& \frac{1}{2} t_d^{-3/2}+\frac{1}{2} t_d^{-3/2}=t_{d+1}^{-1}
\end{eqnarray*}

Lemma \ref{lem:NM1} is proved.\QED

{\it Proof of Lemma \ref{lem:NM2} : } As in the proof of the
previous lemma, the letter $M_k$ is a positive constant which does
not depend on $d$ and which varies from inequality to inequality.\\

Proof of (i) : In the same way as in the proof of the point
$(1_0)$ of the previous lemma, we can show that for all $d\geq
d_k$, we have
\begin{eqnarray*}
\|{\hat \chi}^{d+1}\|_{k+1+s,r_d} &\leq& M_k
\|\zeta(f^d)\|_{k+1+2s,r_d}P(\|f^d\|_{k,r_d})\\
   &+& M_k \|f^d\|_{k+1+2s,r_d}\|\zeta(f^d)\|_{k,r_d}
Q(\|f^d\|_{k,r_d})\\
   &\leq& M_k
\|\zeta(f^d)\|_{k,r_d}^{\frac{k-2s-2}{k-1}}
\|\zeta(f^d)\|_{2k-1,r_d}^{\frac{2s+1}{k-1}} \\
   &+& M_k
\|f^d\|_{k,r_d}^{\frac{k-2s-2}{k-1}}\|f^d\|_{2k-1,r_d}^{\frac{2s+1}{k-1}}
\|\zeta(f^d)\|_{k,r_d}\\
   &\leq& M_k (t_d^{-\frac{k-2s-2}{k-1}+A\frac{2s+1}{k-1}}
+ t_d^{-1+A\frac{2s+1}{k-1}})\\
   &\leq& M_k t_d^{-\mu}
\end{eqnarray*}

where $-\mu<-1/2$. Thus, there exists $d_{k+1} > d_k$ such that
for all $d\geq d_{k+1}$ we have $\|{\hat
\chi}^{d+1}\|_{k+1+s,r_d}<t_d^{-1/2}$. Note that we also have
$\|\chi^{d+1}\|_{k+1+s,r_d}<t_d^{-1/2}$.\\

Proof of (ii) : For $d\geq d_{k+1}$, we have by
(\ref{axiom:action3})
$$
\|f^{d+1}\|_{k+1,r_{d+1}}\leq
\|f^d\|_{k+1,r_d}(1+P(\|{\hat \chi}^{d+1}\|_{k+1+s,r_d}))
$$
where $P$ is a polynomial with vanishing constant term. In Point
(i) we saw that $\|{\hat \chi}^{d+1}\|_{k+1+s,r_d}<t_d^{-1/2}$
then, we can assume, replacing $d_{k+1}$ by a larger integer if
necessary, that $P(\|{\hat \chi}^{d+1}\|_{k+1+s,r_d})\leq
\frac{1}{(d+1)(d+3)}$. Now we choose a positive constant ${\tilde
C}_{k+1}$ (independent on $d$) such that
$\|f^{d_{k+1}}\|_{k+1,r_{d_{k+1}}}<{\tilde
C}_{k+1}\frac{d_{k+1}+1}{d_{k+1}+2}$. We then obtain, as in the
proof of Point $(2)$ of the previous lemma, that
$\|f^{d}\|_{k+1,r_{d+1}}< {\tilde C}_{k+1}\frac{d+1}{d+2}$ for any
$d\geq d_{k+1}$. Note that ${\tilde C}_{k+1}$ is a priori not the
constant of statement of the lemma. Later in the proof (see the
proof of the point (iii)), we will replace it by a
larger one.\\

Proof of (v) : The proof follows the same idea as the proof of
Point $(5)$ in the previous lemma. Let $d$ be an integer such that
$d\geq d_{k+1}-1\geq d_k$.

We have
$$
\|\zeta(f^{d+1})\|_{k+1,r_{d+1}} \leq \|\zeta({\hat
\Phi}_{d+1}\cdot f^{d}-\Phi_{d+1}\cdot f^{d})\|_{k+1,r_{d+1}}
+\|\zeta(\Phi_{d+1}\cdot f^{d})\|_{k+1,r_{d+1}}
$$
Writing (\ref{eqn:estimate-zeta}) with Point (i) and the estimate
$\|f^d\|_{k+1,r_d}<{\tilde C}_{k+1}$, and the interpolation
inequality (\ref{eqn:interpolNM}), we get
\begin{eqnarray*}
\|\zeta(\Phi_{d+1}\cdot f^{d})\|_{k+1,r_{d+1}} &\leq&
\|\zeta(f^d)\|_{k+1+s,r_d}^2 
 R_k(\|f^d\|_{k+s+1,r_d},\|\chi^{d+1}\|_{k+1+s,r_d},\|f^d\|_{k+1,r_d})\\
     &\leq& M_k\|\zeta(f^d)\|_{k,r_d}^{2\frac{k-s-2}{k-1}}
\|\zeta(f^d)\|_{2k-1,r_d}^{2\frac{s+1}{k-1}} \times\\
  & & R_k(\|f^d\|_{k,r_d}^{\frac{k-s-2}{k-1}}
\|f^d\|_{2k-1,r_d}^{\frac{s+1}{k-1}},\|\chi^{d+1}\|_{k+1+s,r_d},\|f^d\|_{k+1,r_d})\\
 &\leq& M_k 
t_d^{-2\frac{k-s-2}{k-1}+2A\frac{s+1}{k-1}+A\tau \frac{s+1}{k-1}}
\end{eqnarray*}
($\tau$ is the degree in the first variable of the polynomials
$R_k$ introduced in Theorem \ref{thm:Nash-Moser}).
Then, by (\ref{cond:lsespsilon}) and (\ref{cond:Aepsilon}), we
have $\|\zeta(\Phi_{d+1}\cdot f^{d})\|_{k+1,r_{d+1}}\leq M_k
t_d^{-\mu}$ where $-\mu<-3/2$ and, replacing $d_{k+1}$ by a larger
integer if necessary, we have $\|\zeta(\Phi_{d+1}\cdot
f^{d})\|_{k+1,r_{d+1}} < \frac{1}{2}t_d^{-3/2}$.\\

On the other hand, exactly in the same way as in the previous
lemma (using the interpolation inequality with $k$ and $2k-1$), we
can show that
\begin{eqnarray*}
\|{\hat \Phi}_{d+1}\cdot f^{d}-\Phi_{d+1}\cdot
f^{d}\|_{k+1,r_{d+1}} &\leq& M_k \|f^d\|_{k+1+s,r_d} \|{\hat
\chi}^{d+1}-\chi^{d+1}\|_{k+1+s,r_d}\\
   &\leq& M_k \|f^d\|_{2k-1,r_d}^{\frac{s+1}{k-1}} \|{\hat
\chi}^{d+1}-\chi^{d+1}\|_{k+1+s,r_d}\\
   &\leq& M_k t_d^{A\frac{s+1}{k-1}} t_d^{-1}
( t_d^{-\frac{k-2s-3}{k-1}+A\frac{2s+2}{k-1}}+t_d^{-1+A\frac{2s+2}{k-1}})\\
   &<& M_k t_d^{-\mu}
\end{eqnarray*}
with $-\mu<-3/2$. Then, using (\ref{eqn:proj}) and replacing $d_{k+1}$ by a larger integer if
necessary, we can write
$$
\|\zeta({\hat \Phi}_{d+1}\cdot f^{d}-\Phi_{d+1}\cdot
f^{d})\|_{k+1,r_{d+1}} < \frac{1}{2} t_d^{-3/2}\,.
$$
We then obtain for all $d\geq d_{k+1}-1$,
$$
\|\zeta(f^{d+1})\|_{k+1,r_{d+1}}<t_{d+1}^{-1}\,.
$$

Proof of (iii) and (iv) : We first write, using the inequality (\ref{eqn:proj}), for all $d\geq d_{k+1}$,
$\|\zeta(f^d)\|_{2k+1,r_d} \leq \|f^d\|_{2k+1,r_d} T_{2k+1}(\|f^d\|_{k+1,r_d})$
where $T_{k+1}$ is a polynomial. Putting $V_{k+1}:=\max (1,T_{k+1}({\tilde C}_{k+1}))$,
we obtain by Point (ii),
\begin{equation}
\|\zeta(f^d)\|_{2k+1,r_d} \leq V_{k+1}\|f^d\|_{2k+1,r_d}\,.
\label{eqn:zetaf2k+1}
\end{equation}
We will use this inequality at the end of the proof.\\

In the same way as in the proof of $(3_d)$ of the
previous lemma, we can show that for all $d\geq d_{k+1}$ we have
\begin{eqnarray*}
\|f^{d+1}\|_{2k+1,r_{d+1}} &\leq& M_k (\|f^d\|_{2k+1,r_d}+\|{\hat
\chi}^{d+1}\|_{2k+1+s,r_d} \|f^d\|_{k+1,r_d})\\
     &\leq& M_k (\|f^d\|_{2k+1,r_d}+\|{\hat
\chi}^{d+1}\|_{2k+1+s,r_d}) \quad {\mbox { by (ii)}}\,.
\end{eqnarray*}
 By (\ref{axiom:smoothing1}) and (\ref{eqn:estimate-chif}), we
write
\begin{eqnarray*}
\|{\hat \chi}^{d+1}\|_{2k+1+s,r_d} &\leq& M_k t_d^{3s+2} \|\chi^{d+1}\|_{2k-1-2s,r_d}\\
   &\leq& M_k t_d^{3s+2} \|\zeta(f^d)\|_{2k-1-s,r_d} P(\|f^d\|_{k,r_d})\\
   & & + \|f^d\|_{2k-1-s,r_d}\|\zeta(f^d)\|_{k,r_d}
   Q(\|f^d\|_{k,r_d})\\
   &\leq& M_k t_d^{3s+2} (\|\zeta(f^d)\|_{2k-1,r_d}
   + \|f^d\|_{2k-1,r_d}\|\zeta(f^d)\|_{k,r_d} )\\
   &\leq& M_k t_d^{A+3s+2}\,.
\end{eqnarray*}
Now, since $A=6s+5>6s+4$, replacing $d_{k+1}$ by a larger integer
if necessary,  we can assume that for any $d\geq d_{k+1}$, we have
$M_k t_d^{A+3s+2}<\frac{1}{2V_{k+1}} t_d^{3A/2}$ (note that it also
implies $M_k< \frac{1}{2V_{k+1}} t_d^{A/2}$). This gives
\begin{equation}
\|f^{d+1}\|_{2k+1,r_{d+1}} \leq \frac{1}{2V_{k+1}}
t_d^{A/2}\|f^d\|_{2k+1,r_d} +\frac{1}{2V_{k+1}} t_d^{3A/2}\,.
\label{eqn:recf}
\end{equation}

We choose a positive constant $C_{k+1}$ such that
$$
C_{k+1}>\max \Big(1,{\tilde C}_{k+1},
\frac{\|f^{d_{k+1}}\|_{2k+1,r_{d_{k+1}}}}{t_{d_{k+1}}^A}\Big)\,.
$$
We then have
$\|f^{d_{k+1}}\|_{2k+1,r_{d_{k+1}}}<C_{k+1}{t_{d_{k+1}}^A}$ and,
using (\ref{eqn:recf}) we obtain by induction :
$$
\|f^d\|_{2k+1,r_d} < \frac{C_{k+1}}{V_{k+1}} t_d^A < C_{k+1}t_d^A\,,
$$
for all $d\geq d_{k+1}$.

Now, by (\ref{eqn:zetaf2k+1}), we have
$$
\|\zeta(f^d)\|_{2k+1,r_d} \leq V_{k+1}\frac{C_{k+1}}{V_{k+1}} t_d^A\,,
$$
for all $d\geq d_{k+1}$.

Moreover, the definition of $C_{k+1}$ completes the proof of the
point
(i), (ii) and (v).\\

Lemma \ref{lem:NM2} is proved. \QED

{\bf Acknowledgements}. We would like to thank the referee for his
careful reading of this paper, and for his numerous remarks which
helped us improve the paper significantly. In particular, the
appendix was added after he challenged us to find a ``more
conceptual'' proof of our main results.

\bibliographystyle{amsalpha}

\providecommand{\bysame}{\leavevmode\hbox
to3em{\hrulefill}\thinspace}

\end{document}